\documentclass[a4paper]{article}
\usepackage{geometry}
\usepackage{amsmath,amssymb}
\usepackage{multicol}
\usepackage{commath}
\usepackage{subfig}
\usepackage{authblk}
\usepackage{todonotes}
\usepackage{hyperref}
\usepackage{booktabs}
\usepackage{doi}
\usepackage[square,numbers]{natbib}
\hypersetup{pdfborder = 0 0 0, colorlinks=false,
 linkcolor=black,citecolor=black, filecolor=black, urlcolor=black, }

\usepackage{algorithm,algpseudocode}
\usepackage[title]{appendix}

\usepackage{charter}
\usepackage{xcolor}

\definecolor{dodgerblue}{rgb}{0.12, 0.56, 1.0}
\definecolor{deepsaffron}{rgb}{1.0, 0.6, 0.2}
\definecolor{pastelviolet}{rgb}{0.8, 0.6, 0.79}
\definecolor{orchid}{rgb}{0.85, 0.44, 0.84}

\newcommand{\revOne}[1]{#1}
\newcommand{\revTwo}[1]{#1}
\newcommand{\revMine}[1]{#1}

\usepackage{cleveref} 

\usepackage{bm}

\newcommand{\pars}[1]{{\left( #1 \right)}}

\newcommand{\nhat}{{\hat n}}
\newcommand{\bdry}{{\partial\Omega}}
\newcommand{\domain}{\Omega}
\newcommand{\hmax}{{h_{\text{max}}}}
\newcommand{\Reps}{{R_\varepsilon}}
\newcommand{\Ceps}{{C_\varepsilon}}
\newcommand{\im}{{\mathrm{i}}}
\renewcommand{\Re}{\operatorname{Re}}
\renewcommand{\Im}{\operatorname{Im}}
\newcommand*{\mytag}{\addtocounter{equation}{1}(\arabic{equation})\quad}
\newcommand{\genvar}{\rho}

\providecommand{\keywords}[1]
{
  \small	
  \textbf{\textit{Keywords---}} #1
}

\title{An adaptive kernel-split quadrature method for parameter-dependent layer potentials}

\author[1]{Fredrik Fryklund\thanks{Corresponding author. E-mail address: \href{mailto:ffry@kth.se}{ffry@kth.se}}}
\author[1]{Ludvig af Klinteberg}
\author[1]{Anna-Karin Tornberg}

\affil[1]{Department of Mathematics, KTH Royal Institute of Technology, Stockholm, Sweden}

\date{}

\begin{document}

\maketitle

\begin{abstract}
  Panel-based, kernel-split quadrature is currently one of the most
  efficient methods available for accurate evaluation of singular and
  nearly singular layer potentials in two dimensions. However, it can
  fail completely for the layer potentials belonging to the modified
  Helmholtz, modified biharmonic and modified Stokes equations. These equations depend
  on a parameter, denoted $\alpha$, and kernel-split quadrature loses
  its accuracy rapidly when this parameter grows beyond a certain
  threshold.
  This paper describes an algorithm that remedies this
  problem, using per-target adaptive sampling of the source
  geometry. The refinement is carried out through recursive bisection,
  with a carefully selected rule set.
  This maintains accuracy for a wide range of the parameter $\alpha$, at an
  increased cost that scales as $\log\alpha$.
  Using this algorithm allows kernel-split quadrature to be both
  accurate and efficient for a much wider range of problems
  than previously possible.\\\\
\end{abstract}

\section{Introduction}

This paper  presents an extension of the panel-based, kernel-split
quadrature scheme by \citeauthor{Helsing2008} \cite{Helsing2008}, for
evaluating singular and nearly singular layer potentials in two
dimensions. It is one of the current state of the art
methods for maintaining low errors when solving homogeneous elliptic partial
differential equations (PDEs) in two dimensions using integral
equation methods \cite{Helsing2008,Helsing2015,Ojala2015}.
However, there exists a set of problems for which this scheme can fail
completely. This includes the following PDEs in $\mathbb{R}^{2}$:
\begin{align}
  (\Delta - \alpha^2)u = 0, & \quad \text{modified Helmholtz},\label{eq:modhelm} \\
  \Delta(\Delta - \alpha^2)u = 0, & \quad \text{modified biharmonic}, \label{eq:modbi}\\
  (\Delta - \alpha^2)u - \nabla p = 0, & \quad \text{modified Stokes (subject to $\nabla \cdot u = 0$)},\label{eq:modstokes}
\end{align}
where $\alpha$ is a positive real number. For brevity, we refer to them as the
\emph{modified} PDEs. Note that they are not consistently named in the
literature. For example, the modified Helmholtz equation is also known
as the screened Poisson equation, the Yukawa equation, the linearized
Poisson-Boltzmann equation, and the Debye-H\"{u}ckel
equation. Meanwhile, the modified Stokes equations are also known as
the Brinkman equations. These PDEs appear in many different
applications: electrostatic interactions in protein and related
biological functions, macroscopic electrostatics and fluid flow on the
micro scale, to mention a few
\cite{Zhou2018,Vorobjev2019,Jiang2013,KROPINSKI2011425,chen_jiang_chen_yao_2015,Helsing2018}.

A common trait of the modified PDEs is that their associated layer
potentials have kernels that either decay exponentially, or have
components that decay exponentially, with a rate that is proportional
to $\alpha$. This decay presents a problem for the above mentioned
kernel-split quadrature. In short, the quadrature method is based on
writing the kernel on a form with smooth functions multiplying
explicit singularities, and then integrating each term separately. In
order to be accurate, these smooth functions have to be locally well
approximated by polynomials, which is increasingly difficult for larger values of $\alpha$.

Large values of $\alpha$ are of interest in the context of elliptic marching. In elliptic marching a semi-implicit temporal discretization is applied to the governing equations. A time-step then involves solving a sequence of elliptic equations, such as the modified PDEs:  e.g. the 
heat equation, the time-dependent
Stokes and Navier-Stokes equations correspond to the modified Helmholtz equation, the modified biharmonic equation and the modified Stokes equations, respectively \cite{kropinskiheateq,fryklundheat,Greengard1998a,AFKLINTEBERG2020109353}. Regardless of the specifics of the time-discretization scheme the \revMine{resulting} equations involve the parameter $\alpha$, where $\alpha^2$ is inversely proportional to the time-step size. The smaller the time-step the lower temporal error, but the spatial problem becomes harder to solve accurately. This is a \revTwo{new} problem in the context of elliptic marching, combined with boundary integral methods applied to the resulting modified PDEs. Earlier, high order and accurate schemes were not available, thus small time-steps were redundant \revTwo{since the spatial error dominated, and the low-error regimes were not feasible to explore on standard desktop computers.} Today efficient methods are available that give
higher precision and allow for more complicated problems. It is
 evident that the restriction of the time-step is a significant
bottleneck, as the kernel-split quadrature may fail for high temporal resolutions.

We have developed a robust quadrature scheme, based on adaptive
refinement, that maintains high accuracy for any $\alpha$, without
sacrificing efficiency. It applies to target points both on and close to the boundary, where
regular quadrature is insufficient. In this context, refinement refers
simply to an interpolation of known quantities to a locally refined
discretization, as opposed to increasing the number of degrees of
freedom in the discretized integral equation. The additional cost, in terms of assembly time per panel, scales as $\mathcal{O}(\log \alpha)$.

The remainder of this paper is organized as follows. In section
\ref{s:background} we give an outline of the kernel-split
quadrature. Section \ref{s:problem} describes the problem that we are
trying to solve, using the modified Helmholtz equation as an
example. In section \ref{s:errorintro} we present error analysis, again with the modified Helmholtz equation as template. The new algorithm we propose is presented in section
\ref{s:algorithm}, followed by numerical results in section
\ref{s:results}.

\section{Background}
\label{s:background}
Our goal is to evaluate layer potentials in $\mathbb{R}^{2}$ of the form
\begin{align}
  u(x) = \int_\bdry G(x, y) \sigma(y) \dif S(y) .
  \label{eq:laypot}
\end{align}
The layer density $\sigma(y)$ is assumed to be smooth, \revOne{implying that the boundary is assumed to be smooth as well}. The
kernel $G(x,y)$ is singular at $x=y$ and can be expressed with
explicit singularities as
\begin{align}
  G(x, y) = G^S(x, y) + G^L(x, y) \log\abs{y-x}
  + G^C(x, y) \frac{(y-x)\cdot\nhat(y)}{\abs{y-x}^2},
  \label{eq:split}
\end{align}
where $G^S$ is a smooth function. The functions $G^L$ and $G^C$ are smooth functions that multiply a log-type singularity and a Cauchy-type singularity, respectively. We refer to this decomposition as \emph{kernel-split}. 

The boundary $\bdry$ is discretized using a composite Gauss\revTwo{--}Legendre
quadrature. It is subdivided into intervals $\Gamma_{i}$, denoted
\emph{panels},
\begin{align}
  \bdry = \bigcup_i \Gamma_i.
\end{align}
Each panel $\Gamma_i$ is described by a parametrization
$\gamma_i$,
\begin{align}
  \Gamma_i = \left\{\gamma_i(t) \in \mathbb{R}^2\vert \, \quad t \in [-1,1]\right\}.
\end{align}
We refer to $x$  in (\ref{eq:laypot}) as a \emph{target point} and a
point $y$ as a \emph{source point}.  A panel to which a source point
belongs is referred to as a \emph{source panel}. Associated with the parametrization is a speed function
$s_i(t) = \abs{\gamma_i'(t)}$, a normal vector $\nhat_i(t)$ and the curvature $\kappa_i(t)$.
Introducing the convenience notation $\sigma_i(t) = \sigma(\gamma_i(t))$,
the layer potential from a panel $\Gamma_i$ becomes
\begin{align}
    \int_{\Gamma_i} G(x, y) \sigma(y) \dif S(y)
    &= \int_{-1}^1 G\pars{x, \gamma_i(t)} \sigma_i(t) s_i(t) \dif t .
\end{align}

Each panel is discretized in the parametrization variable $t$ using
the nodes and weights $(t_j^G, \lambda_j^G)$ of an $n$-point
Gauss\revTwo{--}Legendre quadrature rule on the canonical interval $[-1,\,1]$, which is of order $2n$, such that on
each panel we have the discrete quantities
\begin{align*}
\mytag  y_{ij} &= \gamma_i(t_j^G), &
\mytag  \nhat_{ij} &= \nhat_i(t_j^G), &
\mytag  \sigma_{ij} &= \sigma\pars{\gamma_i(t_j^G)}, \\
\mytag  s_{ij} &= \abs{\gamma_i'(t_j^G)},&
\mytag  \kappa_{ij} &= \kappa_{i}(t_j^G).
\end{align*}
Omitting the panel index $i$, the layer potential contribution from a
panel $\Gamma$ is then computed using the approximation
\begin{align}
  \int_{\Gamma} G(x, y) \sigma(y) \dif S(y)
  &\approx \sum_{j=1}^n G\pars{x, y_j} \sigma_j s_j \lambda_j^G.
    \label{eq:direct_quad}
\end{align}
Due to the singularities in $G$, the above formula requires $x$ to be
\emph{well-separated} from $\Gamma$ in order to be accurate (see section \ref{s:algorithm}, paragraph Near evaluation criterion). Otherwise, the scheme of \cite{Helsing2015} is used, known
as \emph{product integration}. With that, target-specific quadrature weights $w^{L}$ and $w^{C}$ of order $n$ are computed for the known singularities where needed,
such that
\begin{align}
  \int_\Gamma f(x,y) \log\abs{y-x} \dif S(y)
  &\approx \sum_{j=1}^n f(x, y_j) w_j^L(x),
        \label{eq:prodqlog} \\
  \int_\Gamma f(x,y) \frac{(y-x)\cdot\nhat(y)}{\abs{y-x}^2} \dif S(y)
    &\approx \sum_{j=1}^n f(x, y_j) w_j^C(x) .
      \label{eq:prodqcauchy}
\end{align}
Substituting \eqref{eq:split} into \eqref{eq:laypot} and applying the
above product integration gives the so called \emph{kernel-split
  quadrature scheme}. Depending on the location of
the target point $x$ relative to the source panel $\Gamma$, the
evaluation can be divided up into three different cases:

\begin{enumerate}
\item \textbf{Singular, with self-interaction.} If $x\in\Gamma$ is one
  of the quadrature nodes, $x=y_i$, then the term multiplying $G^C$ is
  smooth, with the limit
  \begin{align}
    \label{eq:limitcurvature}
    \lim_{\substack{x \to y\\x, y \in \bdry}} \frac{(y-x)\cdot\nhat(y)}{\abs{y-x}^2} = -\frac{\kappa(y)}{2},
  \end{align}
  where $\kappa(y)$ is the curvature of $\partial\Omega$ at
  $y$. Applying product integration to the $G^{L}$ term, we get
  \begin{align}
    \begin{split}
      \int_{\Gamma} G(y_i, y) \sigma(y) \dif S(y) \approx
      & \sum_{\substack{j=1\\j\ne i}}^n \left[ G(y_i, y_j) \sigma_j s_j \lambda_j^G
        +G^L(y_i, y_j) \sigma_j \pars{w^L_j(y_{i}) - s_j\log\abs{y_j - y_i} \lambda_j^G} \right]\\
      & +G^S(y_i, y_i)
      + G^L(y_i, y_i) \sigma_i w^L_i(y_{i})
      - G^C(y_i, y_i)\frac{\kappa(y_i)}{2}  .
    \end{split}
        \label{eq:specquad_on}
  \end{align}
  Note that this is the only case where $G^{S}$ is evaluated explicitly, and only in one point. 
\item \textbf{Singular, without self-interaction.} If $x \in \bdry$ is
  either on the source panel $\Gamma$ but not a quadrature node, or on
  a neighboring panel, then the $G^C$ term is still smooth, and we do
  not have to take the limit at $x \to y$ into account. This lets us simplify the above to
  \begin{align}
    \begin{split}
      \int_{\Gamma} G(x, y) \sigma(y) \dif S(y) \approx
      & \sum_{j=1}^n \left[ G(x, y_j) \sigma_j s_j \lambda_j^G
        +G^L(x, y_j) \sigma_j \pars{w^L_j(x) - s_j\log\abs{y_j-x}\lambda_j^G} \right] .
    \end{split}
        \label{eq:specquad_on_nb}
  \end{align}

\item \textbf{Nearly singular case.} If $x$ is close to $\Gamma$, but
  not on a neighboring panel, then we need to address both
  singularities in \eqref{eq:split}. This case occurs when either
  $x \in \domain$ is close to $\bdry$, or when $x \in \bdry$ is on a
  section of the boundary that is distant in arc length or
  disjoint. The layer potential is then evaluated as
  \begin{align}
    \begin{split}
      \int_{\Gamma} G(x, y) \sigma(y) \dif S(y) \approx
      \sum_{j=1}^n
        & \Bigg[ G(x, y_j) \sigma_j s_j \lambda_j^G \\
         + &G^L(x, y_j) \sigma_j \pars{w^L_j(x) - s_j\log\abs{y_j-x}\lambda_j^G} \\
         + &G^C(x, y_j) \sigma_j \pars{w^C_j(x) - s_j  \frac{(y_j-x)\cdot\nhat(y_j)}{\abs{y_j-x}^2}\lambda_j^G} \Bigg].
    \end{split}
        \label{eq:specquad_near}
  \end{align}
\end{enumerate}
Given a target point $x$, the panels on $\bdry$ are partitioned into
two sets: \emph{far} panels, that can be evaluated directly using
\eqref{eq:direct_quad}, and \emph{near} panels, that must be evaluated
(or corrected), using either \eqref{eq:specquad_on},
\eqref{eq:specquad_on_nb}, or \eqref{eq:specquad_near}.

\section{Problem statement}
\label{s:problem}

The kernel-split quadrature scheme outlined above is both efficient
and accurate when applied to the single and double layer potentials of
several PDEs, such as the Laplace, Helmholtz and Stokes equations
\cite{Helsing2008,Helsing2015,Ojala2015} \revOne{on geometries with a smooth boundary}. However, for the modified PDEs \eqref{eq:modhelm}, \eqref{eq:modbi}, and \eqref{eq:modstokes}, the scheme can
fail completely.  All of these equations have layer potential kernels
including second-kind modified Bessel functions, of the forms
$K_0(\alpha\revTwo{|y-x|})$ and/or $K_1(\alpha\revTwo{|y-x|})$. As we shall see, this is
problematic for the kernel-split quadrature.
\begin{figure}[htbp]
  \centering
  \includegraphics[width=0.8\textwidth]{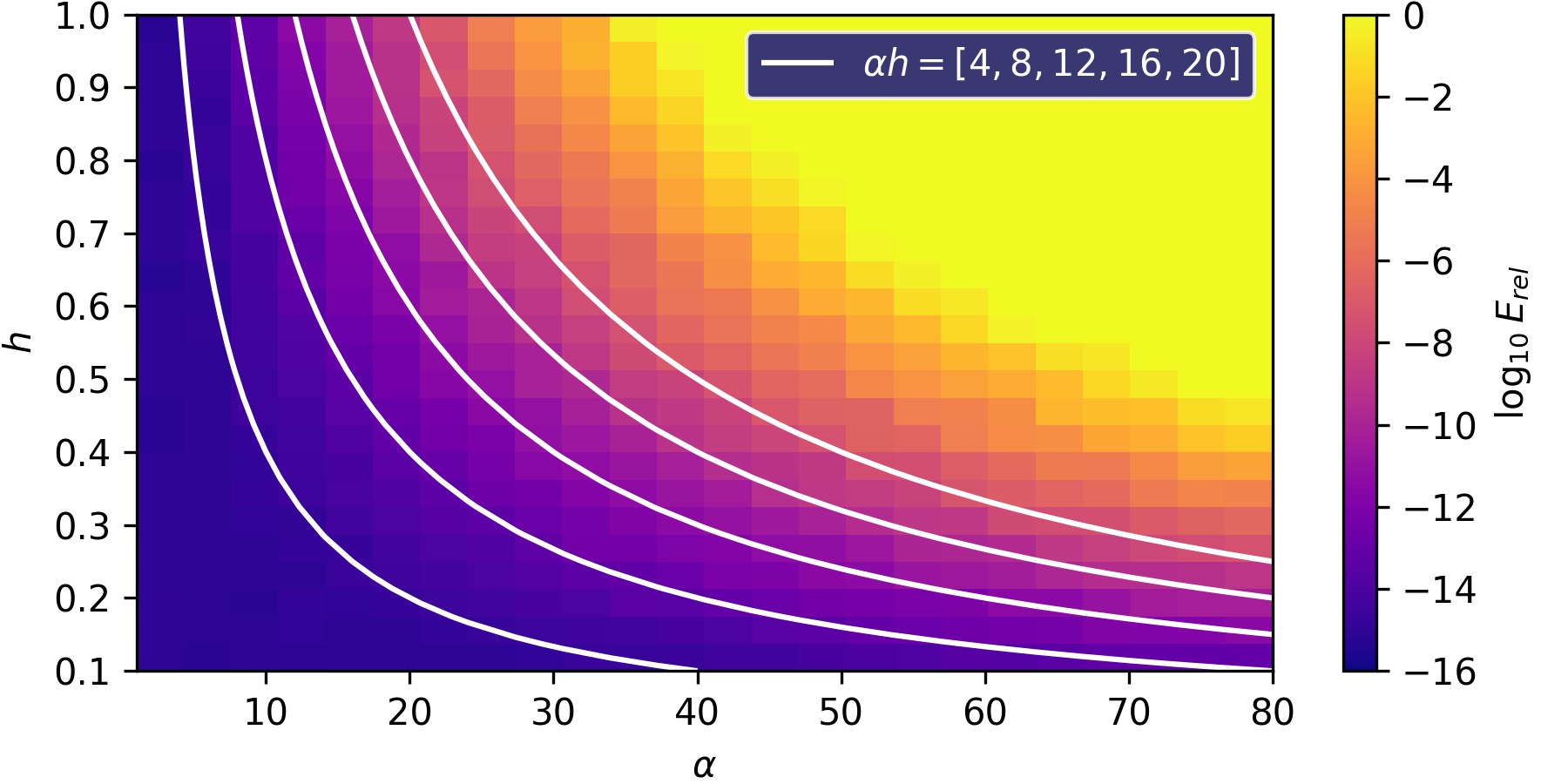}
  \caption{Relative error
    $E_{\mathrm{rel}} = \|u-\tilde u\|_{\infty}/\|u\|_{\infty}$, where $\tilde{u}$ is the result of  evaluating the
    modified Helmholtz single layer potential \eqref{eq:sgl_demo_int}
    over a flat panel of length $h$, using kernel-split
    quadrature with $n = 32$ Gauss\revTwo{--}Legendre points. The norm is taken over 100 values of target points $x$ randomly
    drawn from the box $[-h/2, h/2]\times[0,h/2]$. The white lines are
    contours of the quantity $\alpha h$, providing a heuristic
    motivation for why a criterion of the form \eqref{eq:alpha_h_crit}
    is suitable.}
  \label{fig:ah_demo}
\end{figure}

As an illustrating example, we study the single layer kernel of the
modified Helmholtz equation~\eqref{eq:modhelm}. The associated Green's function is
\begin{align}
  G(x,y) = K_0(\alpha\revTwo{|y-x|}),
  \label{eq:modhelm_green}
\end{align}
where the scaling factor $1/2\pi$ has been omitted. The split of this kernel is based on using a standard decomposition
\cite[\S10.31]{NIST:DLMF} to explicitly write out the singularities in
$K_0$,
\begin{align}
  K_0(\genvar) = K_0^S(\genvar) - I_0(\genvar) \log \genvar,\quad \genvar\in\mathbb{R}^{+},
  \label{eq:K0_split}
\end{align}
where $K_0^S$ is the smooth remainder, and $I_0$ is the modified Bessel function of the first kind. For future reference, we have by \cite[\S10.25,\S10.29]{NIST:DLMF} that
\begin{equation}
\label{eq:defI}
  I_{\nu}(\rho) = \left(\frac{1}{2}\rho\right)^{\nu}\sum\limits_{k=0}^{\infty}\frac{\left(\frac{1}{4}\rho^{2}\right)^{k}}{k!(\nu+k)!},\quad \nu\in\mathbb{N},\,\rho\in\mathbb{R},
\end{equation}
and the derivatives of $I_{\nu}$ can be expressed as
\begin{align}
\label{eq:defIdiff}
  I_{\nu}^{(n)}(\genvar) = 2^{-n} \sum_{m=0}^n \binom{n}{m} I_{2m+(\nu-n)}(\genvar),\quad \nu,\,n\in\mathbb{N},\,\rho\in\mathbb{R},
\end{align}
 which makes them easy to evaluate using standard numerical libraries,

Inserting \eqref{eq:K0_split}  into the kernel
\eqref{eq:modhelm_green}, after also splitting the logarithm, we
identify the terms in \eqref{eq:split} as
\begin{align}
\label{eq:modhelm_G0}  G^S(x,y) &=
             K_0^S(\alpha\revTwo{|y-x|}) \revMine{-} I_0(\alpha \revTwo{|y-x|}) \log \alpha
             , \\
\label{eq:modhelm_GL}  G^L(x,y) &= -I_0(\alpha\revTwo{|y-x|}), \\
\label{eq:modhelm_GC}  G^C(x,y) &= 0 .
\end{align}
Product integration is a semi-analytical method. This means that
the target-specific quadrature weights in \eqref{eq:prodqlog}
are found by an analytic treatment of the
singularity or near singularity. The accuracy of the method for a
given evaluation point $x$ is
basically limited by how well a function $f(x,y)$ can be resolved by an ($n-1$)th-degree Legendre polynomial. For a well-resolved function, the
integration error is
usually only at most a few orders of magnitude larger than round-off,
independent of the location of the evaluation point.
We now want to evaluate the layer potential
(\ref{eq:laypot}) with the layer density $\sigma$, using the split
above. With this $G^L$,
$f(x,y)$ in (\ref{eq:prodqlog}) will be $\revTwo{-}I_0(\alpha\revTwo{|y-x|})
\sigma(y)$. Hence, for good accuracy in integrating the logarithmic
singularity, it is this product that must have a small
interpolation error.

The function $G(x,y)$ goes to zero as $\revTwo{|y-x|}$ goes to infinity,
since $K_{\nu}(\genvar)\sim \sqrt{\pi/(2\genvar)}e^{-\genvar}$ as
$\genvar \to \infty$ for $\nu = 1,2$, but its components $G^S(x,y)$
and $G^L(x,y)$ actually grow exponentially as $e^{\alpha\revTwo{|y-x|}}$, with opposing signs, 
following the asymptotic
$ I_{\nu}(\genvar) \sim e^{\genvar}/\sqrt{2 \pi \genvar}$ as
$\genvar \to \infty$ for $\nu = 1,2$
\cite[\S10.25,\S10.30]{NIST:DLMF}. As $\alpha$ gets larger, this makes
$G^{L}$ an increasingly bad candidate for polynomial interpolation,
which will render the product integration inaccurate. In
addition, when evaluated in limited precision, this split is prone to
numerical cancellation, due to the limiting forms of $K_{0}$ and
$I_{1}$.  We will argue that the error mainly is a function of the
quantity $\alpha h$, which can be heuristically motivated by figure
\ref{fig:ah_demo}.  To ensure that this error remains below some
tolerance $\epsilon$, we therefore suggest that $\alpha$ and \revMine{panel length} $h$ must
satisfy a criterion on the form
\begin{align}
  \alpha h \le \Ceps,
  \label{eq:alpha_h_crit}
\end{align}
for some constant $\Ceps$, which can be determined using numerical
experiments. The above criterion can also be reformulated as follows: In order to
achieve a tolerance $\epsilon$, panel lengths must satisfy
\begin{align}
  h \le \hmax := \Ceps / \alpha . \label{eq:hmax_crit}
\end{align}
The naive way of achieving this, for a given $\alpha$, is to
discretize $\bdry$ using sufficiently short panels. However, this can
result in a discretization with orders of magnitude more points than
necessary to resolve the geometry and the layer density.
 Global refinement is also redundant for another reason: as stated above the functions
 $K_0(\genvar)$ and $K_1(\genvar)$ decay as $\sqrt{\pi/(2 \genvar)}e^{-\genvar}$, meaning that they are very localized for large $\alpha$, i.e. when polynomial interpolation might
 fail, and are almost zero in finite precision away from $\genvar = 0$. Thus only a small portion of the boundary needs refinement.

\section{Error estimates}
\label{s:errorintro}
As was indicated above, there are two main sources of errors in the
kernel split quadrature that grow as $\alpha$ increases. The first
part arises from an interpolation error, and the second is due to
numerical cancellation. In order to better understand these errors, we
will perform a limited analysis for the single layer potential of the
modified Helmholtz equation on a flat panel in section
\ref{ss:erroranalysis}. 

This understanding of the error structure is useful also for the
associated kernels for the other modified PDEs since they have a
similar singularity structure to the single layer modified Helmholtz potential,
see appendix \ref{sec:kernelsplits}. They are all combinations of $K_{0}$ and/or $K_{1}$,
thus in the form of  explicit singularities they contain $I_{0}$ and/or $I_{1}$ multiplying a log-type
singularity. 

Based on our simple analysis in section \ref{ss:erroranalysis}, we will
in section \ref{ss:errorestimate} formulate an error estimate and
illustrate its performance for the single layer
modified Helmholtz potential.
The error estimates can be used to set $C_{\epsilon}$ in 
\eqref{eq:alpha_h_crit}, as discussed in section \ref{ss:Cepsilon}. Other kernels are discussed in section \ref{ss:otherkernels}.

In this section as well as the next, we will identify vectors $x,y \in \mathbb{R}^{2}$
with the corresponding complex numbers $x,y \in \mathbb{C}$.

\subsection{Error analysis for a flat panel}
\label{ss:erroranalysis}

We consider the single layer
potential for the modified Helmholtz equation with unit density, evaluated using kernel-split quadrature
from a flat panel of length $h$,
\begin{align}
  u(x) = \int_{-h/2}^{h/2} K_0(x,y) \dif y.
  \label{eq:sgl_demo_int}
\end{align}
To evaluate \eqref{eq:sgl_demo_int} using the kernel-split correction
\eqref{eq:specquad_near}, we first write
\begin{align}
  u(x) \approx \sum_{j=1}^n \pars{G(x, y_j) - G^L(x, y_j)\log\abs{y_j-x}} \lambda_j
  - \int_{-h/2}^{h/2} I_0(\alpha\revTwo{|y-x|}) \log|y-x| \dif y,
  \label{eq:sgl_demo_prodint}
\end{align}
using $G^C(x,y)= 0$ \eqref{eq:modhelm_GC}. In order to evaluate the remaining integral, the \revTwo{function $I_0(\alpha\revTwo{|y-x|})$} is approximated by an $(n-1)$th-degree polynomial $p_{n-1}$ that interpolates \revTwo{$I_0(\alpha|y-x|)$} at the $n$ nodes $y_j$,
\begin{align}
  I_0(\alpha\revTwo{|y-x|}) &= p_{n-1}(x, y) + r_n(x, y) ,
                     \\
  p_{n-1}(x, y) &= \sum_{k=0}^{n-1} c_k(x) y^{k} \label{eq:poly}, \\
  p_{n-1}(x, y_j) &= I_0(\alpha\revTwo{|y_{j}-x|}), \quad j=1,\ldots,n,
\end{align}
where $r_n$ is the polynomial interpolation error. We consider only even values for $n$. The integral on the right-hand side of \eqref{eq:sgl_demo_prodint} can be written as a sum: the approximation by integrating the polynomial $p_{n-1}$, and an integral over the polynomial interpolation error
\begin{align}
  \int_{-h/2}^{h/2} I_0(\alpha\revTwo{|y-x|}) \log|y-x|\dif y =
                                                        \mathrm{Re}\left( \int_{-h/2}^{h/2} I_0(\alpha\revTwo{|y-x|}) \log(y-x) \dif y\right)\\
  =
  \mathrm{Re}\left(\sum_{k=0}^{n-1} c_k(x) q_k(x)\right) + \mathrm{Re}\left(\int_{-h/2}^{h/2} r_n(x,y) \log(y-x) \dif y\right).
  \label{eq:specquad_w_err}
\end{align}
Here, we have used the definition of the complex logarithm.
We have that
\begin{align}
  q_k(x) = \int_{-h/2}^{h/2} y^k\log(y-x)\dif y, \quad k = 0,\ldots,n-1,
  \label{eq:qk_def}
\end{align}
are integrals that can be computed recursively, starting from exact
formulas \cite{Helsing2009}.

The error in the kernel-split quadrature has two sources. The first is the integral over the interpolation error from the
integral in the right-hand side of \eqref{eq:sgl_demo_prodint}, namely

\begin{align}
  R_n(\alpha,x,h) = \int_{-h/2}^{h/2} r_n(x,y) \log|y-x|\dif y =  \text{Re}\left(\int_{-h/2}^{h/2} r_n(x,y) \log(y-x)\dif y\right).
\label{eq:R_integral}
\end{align}
For brevity we also refer to the integral over the interpolation error as the interpolation error. It will be clear from context which one that is referred to.

The second is a cancellation error due to the limiting forms of $G$ and $G^{L}$. Recall that $G^{S} = G - G^{L}$; this function grows exponentially as $-e^{\alpha\revTwo{|y-x|}}/\sqrt{2\pi\alpha\revTwo{|y-x|}}$, which is the same
asymptotics as for $G^{L}$, but with opposite sign. The function $G$ is a decaying function as $\alpha\revTwo{|y-x|}$ goes to infinity, thus $G^{S}$ and $G^{L}$ have to cancel. As the quadrature terms in \eqref{eq:sgl_demo_prodint} grows in magnitude catastrophic cancellation errors follow.

We will now perform the analysis of the interpolation error for target
points $x$ along the real axis. In section \ref{ss:errorestimate}, we will
consider the performance of this estimate for target points in the
plane.

\subsubsection{Interpolation error}

Our goal is to estimate the magnitude of the interpolation error
$R_n(\alpha,x,h)$ in (\ref{eq:R_integral}), and we will do so for $x
\in \mathbb{R}$.
The polynomial $p_{n-1}$ \eqref{eq:poly} is formed by interpolating at
the Legendre nodes, $y_j=t_j^Gh/2$. According to \cite{PowellM.J.D1981Atam}, the remainder is given by
\begin{equation}
  r_n(x,y) = \prod_{j=1}^n (y-t_j^Gh/2) \frac{1}{n!} \dod[n]{}{\xi} I_0\pars{\alpha(
  \xi-x)}, \quad y,\,\xi\in[-h/2,\,h/2],
\label{eq:standard_interp_error}
\end{equation}
where the absolute value of the argument of $I_{0}$ can be removed, since $I_{0}$ is an even function. This is convenient, as it \revMine{simplifies differentiation}.

The $n$ Legendre nodes are the roots of the Legendre polynomial of order $n$, denoted $P_{n}$, orthogonal on $[-h/2,\,h/2]$. Hence
\begin{align}
  \ell_n \prod_{j=1}^n (y-t_j^Gh/2) = P_n(2y/h),\quad y\in [-h/2,\,h/2],
  \label{eq:legprod}
 \end{align}
 where $\ell_n$ is the coefficient of the leading order monomial term in $P_n$. From Rodrigues' formula \cite[\S18.5]{NIST:DLMF}, it can be shown that $\ell_n=2^{-n}(2n)!/(n!)^2$. After rescaling $\xi$ we have
\begin{align}
  r_n(x,y) = P_n(2y/h)\frac{(\alpha h)^{n}n!}{(2n)!} I_0^{(n)}\pars{\alpha(\xi  h/2 - x)}\, \quad \xi\in[-1,1].
\label{eq:standard_interp_error_bound}
\end{align}
To present an error bound for \eqref{eq:standard_interp_error_bound} we can apply
\begin{equation}
  0\leq I_0^{(n)}\pars{\alpha(\xi h/2 - x)} \leq
  I_0^{(n)}\pars{\alpha(h/2 + |x|)},\quad \forall \xi\in[-1,1].
\label{eq:errbound}
\end{equation}
However, the resulting error bound greatly overestimates the error, thus is not of practical use. We pursue an alternative approach by  estimating $r_{n}$ by choosing a value for $\xi$ for all pairs $(x,y)$. For this purpose, heuristics suggest $\xi = 0$, as demonstrated below.

Combining \eqref{eq:R_integral} and \eqref{eq:standard_interp_error_bound} for $\xi = 0$ we get the estimate
\begin{align}
  \abs{R_n(\alpha,x,h)} \approx \abs{\text{Re}\left(\Lambda_{S}(\alpha, n , x, h)\right)}\leq  \abs{\Lambda_{S}(\alpha, n , x, h)}
                        ,
                        \label{eq:R_est_Lambda}
\end{align}
with $\Lambda_{S}$ defined as
\begin{align}
  \Lambda_{S}(\alpha, n , x, h) &=
                  \frac{n!}{(2n)!}
                  (\alpha h)^n
                  I_0^{(n)}\left(\alpha |x|\right)\frac{h}{2}
g_n\left(2x/h\right),
\label{eq:lambdadef}
\end{align}
where the subscript $S$ denotes single layer potential. Note that $I_{0}^{(n)}$ is an even function for even values of $n$;  by the relation \eqref{eq:defIdiff} $I_{0}^{(n)}$ is a linear combination of $I_{\nu}$ with even $\nu$, and it follows from the definition  \eqref{eq:defI} that these terms are even. Thus we may take the absolute value of the argument of $I_{0}^{(n)}$ without altering the result. This is important, as later on we will consider complex target points $x$.

In \eqref{eq:lambdadef} we have
\begin{align}
  \label{eq:leglogint}
g_{n}(\tilde{x}) = \int_{-1}^{1} P_n(\genvar) \log(\genvar-\tilde{x}) \dif \genvar
\end{align}
obtained by
\begin{align}
\int_{-1}^{1} P_n(\genvar) \log(\genvar h/2-x) \dif \genvar
 = \int_{-1}^{1} P_n(\genvar) \log(\genvar-2x/h) \dif \genvar = g_{n}(2x/h),
\end{align}
where the integral is rewritten using
$\log(h/2(\genvar-2x/h))=\log(h/2)+\log(\genvar-2x/h)$ and the fact
that each $P_n(\rho)$ integrates to zero over the interval.

This ``Legendre-log integral'' can be evaluated using the
recursion formulas in section \ref{sec:leglog}.

\begin{figure}[htbp]
  \centering
  \includegraphics[width=0.49\textwidth]{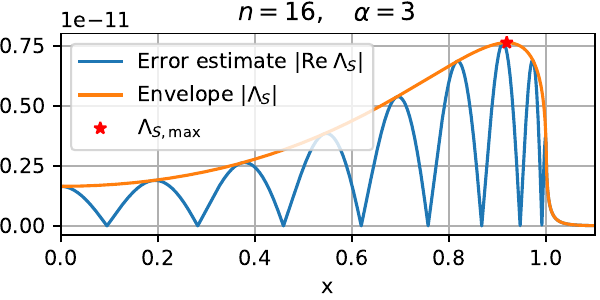}
  \hfill
  \includegraphics[width=0.49\textwidth]{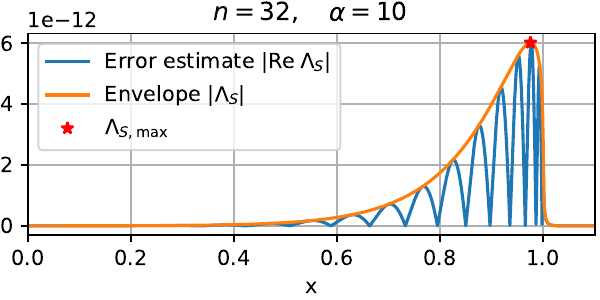}
  \caption{Plot of the error estimate~\eqref{eq:R_est_Lambda} with
    $(n,\alpha)$ valued $(16,3)$ and $(32,10)$, and $h = 2$. The estimate
    $|\Re\Lambda_{S}|$ oscillates in $x$, and is enveloped by the smooth
    upper bound $|\Lambda_{S}|$, which has a maximum to the
    left of $x=1$.}
  \label{fig:est_real_line}
\end{figure}

For $h = 2$, both the absolute value as well as the real part of $\Lambda_{S}$ in
(\ref{eq:lambdadef}) are shown as a function of $x$ in figure \ref{fig:est_real_line}. The absolute
value of $\Lambda_{S}$ is a smooth function that achieves its maximum for
$x$ close to the endpoints of the interval.
This is due to the $\log$-factor in the integrand being singular for
target points $x$ in $[-1,1]$, combined with the growth of the factor
$I_{0}$ towards the edges of the panel.  For target points $x$ outside
the source panel the Legendre-log integral $g_{n}$ decreases rapidly.

This plot is only illustrating the error estimate derived for this
simplifed case. How well it actually estimates the error will be
discussed after we have considered also the second part of the error. 

\subsubsection{Numerical cancellation error}
In addition to the interpolation error just discussed, the
kernel-split \eqref{eq:sgl_demo_prodint} can also suffer from
numerical cancellation when $\alpha$ is large. To see this, note that $G\revMine{(x,y)}-G^{L}\revMine{(x,y)}\revTwo{\log(\alpha\abs{y-x})} = G^{S}\revMine{(x,y)}$, which has the asymptotics \revMine{$-e^{\alpha\abs{y-x}}/\sqrt{2 \pi \alpha\abs{y-x}}$}  as $\alpha\abs{y-x}$ goes to \revMine{infinity}. This is the same asymptotics as for $I_{0}$, but with opposing sign. These two terms have to cancel, since $K_{0}$, i.e. $G$, is a decaying function. The sum of terms with large magnitude, but opposing signs, is prone to numerical cancellation in limited precision.

The cancellation error is straightforward to estimate by
\begin{align}
  \Xi_{S}(\alpha,x,h) = \epsilon_{\mathrm{mach}}h I_0\pars{\alpha d},
                    \quad d = \max_{y\in[-h/2,h/2]}|y-x| = h/2 + |x|,
  \label{eq:canc_err_est}
\end{align}
where $\epsilon_{\mathrm{mach}}$ is the machine epsilon. 
\subsection{Error estimate}
\label{ss:errorestimate}
We now have an error estimate  $|\Lambda_{S}(\alpha,n,x,h)|$ \eqref{eq:lambdadef} for the
interpolation error $R_{n}(\alpha,x,h)$ \eqref{eq:R_integral}, and an estimate $\Xi_{S}(\alpha,x,h)$ \eqref{eq:canc_err_est} for the
cancellation error. Both were derived for target points along the
real axis, but we now evaluate the estimated maximum error over a set $D$ of discrete target
points in the plane, and compare it to the maximum measured actual error. There are two goals: to study the errors' and error estimates' dependence on $\alpha h$, and to construct a combination of the error estimates to predict the total error.

\begin{figure}[htbp]
  \centering
  \includegraphics[width=0.49\textwidth]{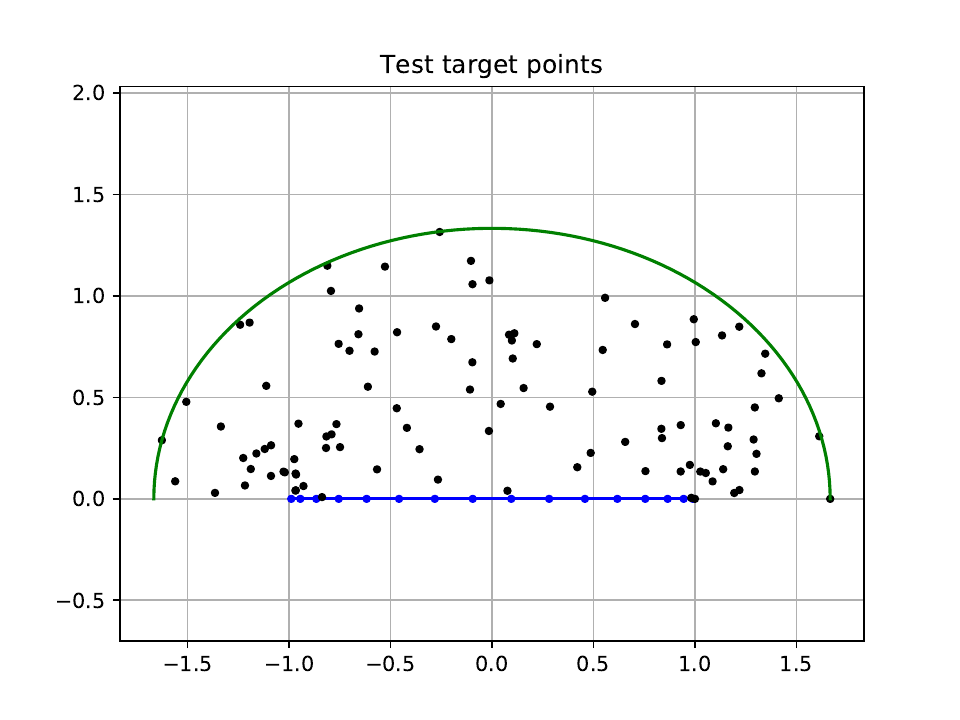}
  \hfill
    \includegraphics[width=0.49\textwidth]{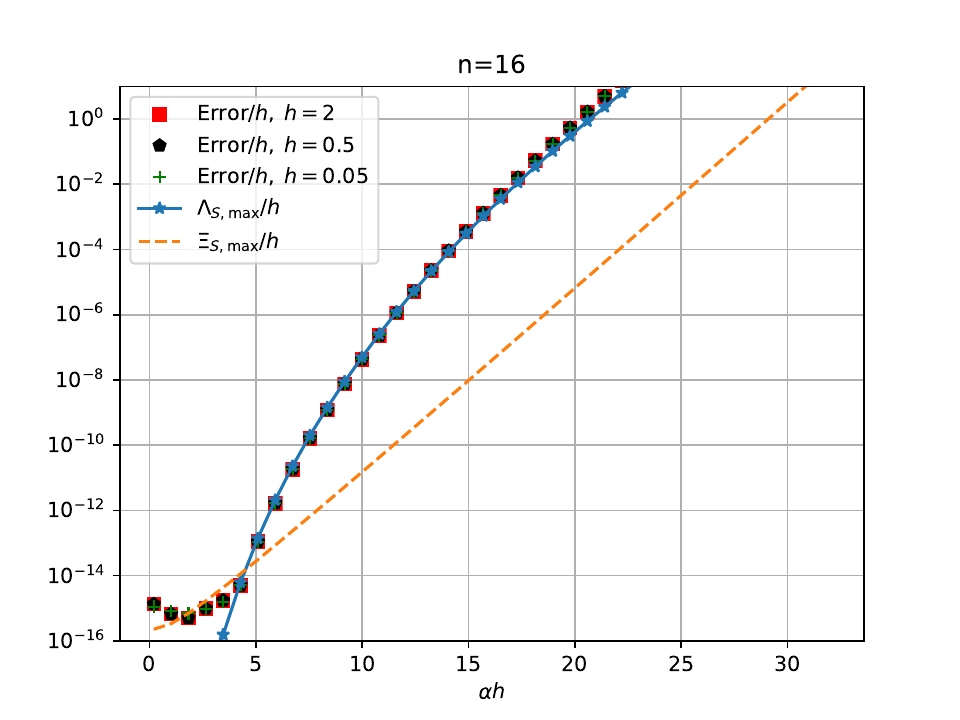}
  \caption{Left: Distribution of $100$ target points in the upper half
    of the Bernstein ellipse with foci at $\pm 1$, forming the set $D$,
    and a flat panel. Right: Plot of measured actual total errors scaled with $1/h$ for different $h$, error bound and error estimates $\Lambda_{S,\max}/h$ and $\Xi_{S,\max}/h$  for the target points in the left plot for the single layer kernel for the modified Helmholtz equation \eqref{eq:modhelm_green}, with
    $n=16$ over a range of $\alpha h$.}
  \label{fig:errestn}
\end{figure}

\begin{figure}[htbp]
  \centering
  \includegraphics[width=0.49\textwidth]{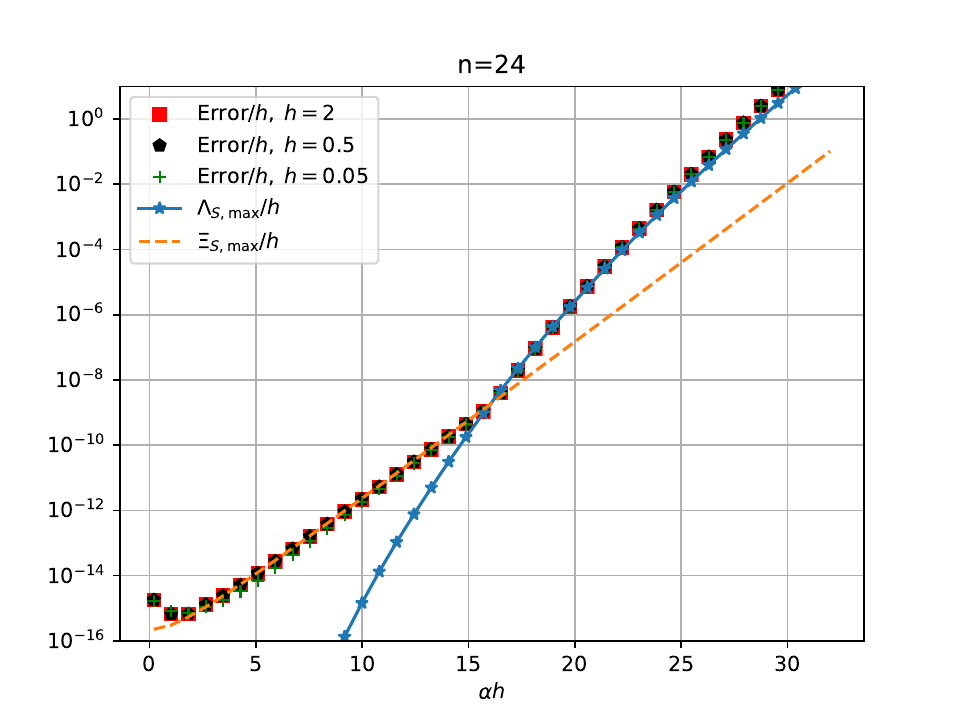}
  \hfill
    \includegraphics[width=0.49\textwidth]{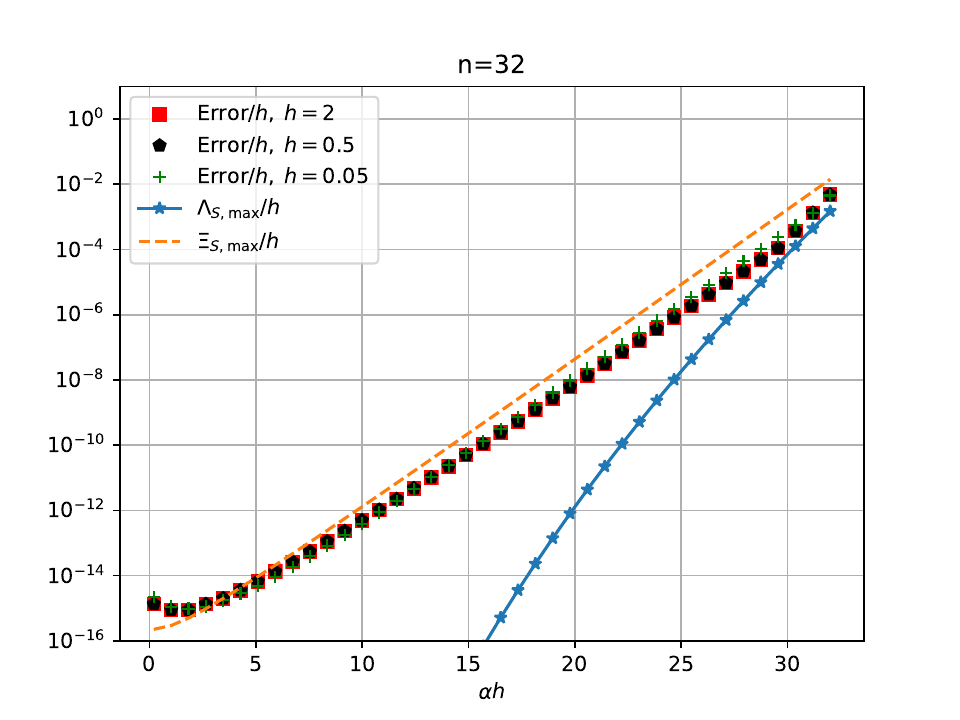}
  \caption{Plots of measured actual total errors scaled with $1/h$ for different $h$, and error estimates $\Lambda_{S,\max}/h$ and $\Xi_{S,\max}/h$  for the target points in figure \ref{fig:errestn} for the single layer kernel for the modified Helmholtz equation \eqref{eq:modhelm_green}, with $n=24$ (left) and $32$ (right), over a range of $\alpha h$.}
  \label{fig:errestn2}
\end{figure}

Let $D$ be the set $100$ sampled complex points with positive imaginary part uniformly within a Bernstein ellipse with foci $\pm 1$, and denote its elements as $\tilde{x}$. The radius for the \revTwo{Bernstein} ellipse is set to be equal to $3^{16/n}$, as points farther away from the panel do not need the kernel-split quadrature scheme to be accurate \cite{klinteberg2020quadrature}. This set is visualized in figure \ref{fig:errestn}. By setting $\tilde{x} = 2x/h$, the estimate \eqref{eq:lambdadef} can be rewritten as 
\begin{align}
  \Lambda_{S}(\alpha, n , \tilde{x}, h) &=
                  \frac{n!}{(2n)!}
                  (\alpha h)^n
                  I_0^{(n)}\left(\frac{\alpha h}{2} |\tilde{x}|\right)\frac{h}{2}g_n(\tilde{x}),
\label{eq:lambdadef2}
\end{align}
and similarly the cancellation error \eqref{eq:canc_err_est} becomes
\begin{align}
  \Xi_{S}(\alpha,\tilde{x},h) = \epsilon_{\mathrm{mach}}h I_0\pars{\frac{\alpha h}{2} d},
                    \quad d = 1+|\tilde{x}|.
\label{eq:canc_err_est_2}
\end{align}
We define 
\begin{align}
\Lambda_{S, \max}&= \max_{\tilde{x}\in D}\, |\Lambda_{S}(\alpha,n,\tilde{x},h)|, 
\label{eq:iperrmax} \\
\Xi_{S,\max}& = \max_{\tilde{x}\in D}\, \Xi_{S}(\alpha,\tilde{x},h).
\label{eq:ximax}
\end{align}
The maximum of these two quantities gives the error estimate
\begin{equation}
  E_{S}(\alpha,n,h) = \max\, \left(\Lambda_{S,\max},\Xi_{S,\max}\right),
  \label{eq:estmax}
\end{equation}
for given values of $\alpha$, $n$ and $h$. From  \eqref{eq:lambdadef2} and \eqref{eq:canc_err_est_2} it clear that $E_{S}$ is a function with  $\alpha h$ and $h$ as two separate arguments, motivating the formulation
\begin{equation}
  E_{S}(\alpha,n,h) = h \tilde{E}_{S}(\alpha h, n).
\label{eq:errorscaled}
\end{equation}
Hence, we expect the error divided by $h$ to be a function of $\alpha h$ only, and this is confirmed in figure \ref{fig:errestn} and figure \ref{fig:errestn2}. Here, we plot the error scaled with $1/h$ for different values of $h$ versus $\alpha h$, for different values of $n$. We see that the error curves collapse for different values of $h$, just as predicted. Furthermore, from the plots it is clear that the interpolation error and cancellation error dominate in different regimes. As $n$ is increased,
the interpolation error decreases, and the cancellation error will
continue to be dominant for larger values of $\alpha h$. 
Considering the left plot of figure \ref{fig:errestn2}, we can most
clearly see how the measured error follows the cancellation error
estimate $\Xi_{S,\max}$ first, and then shifts to follow the interpolation error estimate $\Lambda_{S,\max}$ as it becomes dominant.  We can conclude
that the error estimate $\tilde{E}_{S}$ \revMine{predicts} the actual error\revMine{,} scaled with $1/h$\revMine{,} quite well.

\subsection{Other kernels}
\label{ss:otherkernels}
We now discuss error estimates for the other kernels of interest for the modified PDEs. The double layer kernel for modified Helmholtz is discussed in detail, which is used to motivate a similar treatment of the error estimates for the modified biharmonic equation and the modified Stokes equations. Their corresponding kernels, and the decompositions thereof into explicit singularities, are shown in appendix \ref{sec:kernelsplits}.

The kernels for layer potentials associated with the modified Helmholtz equation (double layer), modified biharmonic and modified Stokes are harder to study analytically than the single layer kernel for the modified Helmholtz equation; $G^{L}$ for all of them contain $I_{0}$ and/or $I_{1}$ multiplying a smooth function, which we collectively denote $F$. This means that the interpolation error \eqref{eq:standard_interp_error} includes the $n$th derivative of a product. However, we observe that the derived error estimates of the interpolation error for the single layer kernel can be applied to the these kernels as well, with a few modifications. This is not surprising, as the difficulties lie in properly resolving $I_{0}$ and/or $I_{1}$, not resolving the function $F$. The cancellation error \eqref{eq:ximax} is also straightforward to modify.

Before proceeding with the error estimates we present the split of the double layer kernel for the modified Helmholtz equation into explicit singularities. The associated kernel is, again omitting the scaling $1/2\pi$,
\begin{align}
  G(x,y) = \frac{\partial}{\partial \hat{n}(x)}K_0(\alpha\revTwo{|y-x|})=-\alpha K_1(\alpha\revTwo{|y-x|})\frac{(y-x)\cdot \hat{n}(y)}{|y-x|},
  \label{eq:modhelm_green_dbl}
\end{align}
and by \cite[\S10.31]{NIST:DLMF} one has
\begin{align}
  K_1(\genvar) = K_1^S(\genvar) + \frac{1}{\genvar}  + I_1(\genvar) \log \genvar, \quad \rho\in\mathbb{R}_{+}.
  \label{eq:K1_split} 
\end{align}
The content of the functions $G^{0}$, $G^{L}$ and $G^{C}$ depend on the location of the target point $x$. If $x$ belongs to the boundary, i.e. $x\in\partial\Omega$, then the kernel-split is
\begin{align}
\label{eq:modhelm_dbl_G0}  G^S(x,y) &=
            -\alpha\left(K_1^S(\alpha\revTwo{|y-x|}) + \frac{1}{\alpha|y-x|}+ I_0(\alpha \revTwo{|y-x|}) \log \alpha\right)\frac{(y-x)\cdot \hat{n}(y)}{|y-x|},\\
\label{eq:modhelm_dbl_GL}  G^L(x,y) &= -\alpha I_1(\alpha\revTwo{|y-x|})\frac{(y-x)\cdot \hat{n}(y)}{|y-x|}, \\
\label{eq:modhelm_dbl_GC}  G^C(x,y) &= 0 .
\end{align}
For off-boundary $x$ the kernel-split is
\begin{align}
\label{eq:modhelm_dbl_G0_domain}  G^S(x,y) &=
            -\alpha \left(K_1^S(\alpha\revTwo{|y-x|}) + I_1(\alpha \revTwo{|y-x|}) \log \alpha\right)\frac{(y-x)\cdot \hat{n}(y)}{|y-x|},\\
\label{eq:modhelm_dbl_GL_domain}  G^L(x,y) &= -\alpha I_1(\alpha\revTwo{|y-x|})\frac{(y-x)\cdot \hat{n}(y)}{|y-x|}, \\
\label{eq:modhelm_dbl_GC_domain}   G^C(x,y) &= -1.
\end{align}
The main difference between the two decompositions is that $1/|y-x|$ in \eqref{eq:K1_split} goes into $G^{S}$ for $x\in\partial\Omega$, due to the known limit value \eqref{eq:limitcurvature}, while for $x\notin\partial\Omega$ it results in a non-zero $G^{C}$. For both kernel-splits the function $G^{L}$ is smooth and consists of the modified Bessel function $I_{1}$ multiplied with $F = (y-x)\cdot \hat{n}(y)/|y-x|$. Here, $F$ does not pose the difficulties with polynomial interpolation as $I_{1}$ does.

We now create an error estimate for the interpolation error. Again identifying the vectors $x$ and $y$ in $\mathbb{R}^{2}$ with complex numbers $x$ and $y$ in $\mathbb{C}$, and following the same steps as for the single layer potential we get analogous to \eqref{eq:standard_interp_error_bound}
\begin{align}
 & r_n(x,y) = P_n(2y/h)\frac{(\alpha h)^{n}n!}{(2n)!} \frac{\mathrm{d}^{n}}{\mathrm{d}\xi^{n}}\left(\alpha I_1\pars{\alpha|\xi  h/2 - x|}\frac{\pars{\xi  h/2 - x}\cdot\hat{n}(\xi h/2)}{|\xi  h/2 - x|}\right)\\
&= P_n(2y/h)\frac{(\alpha h)^{n}n!}{(2n)!} \frac{\mathrm{d}^{n}}{\mathrm{d}\xi^{n}}\left(\mathrm{sgn}\pars{\xi  h/2 - x}\alpha I_1\pars{\alpha\pars{\xi  h/2 - x}}\frac{\pars{\xi  h/2 - x}\cdot\hat{n}(\xi h/2)}{|\xi  h/2 - x|}\right)
  , \quad \xi\in[-1,1],
\label{eq:standard_interp_error_bound_dlp}
\end{align}
since $I_{1}$ is an odd function. The derivation for the single layer kernel that we based this result on was done for $x$ in $\mathbb{R}$. We will now, as before, apply it to target points $x$ in $\mathbb{C}$, for which  $\pars{\xi  h/2 - x}\cdot\hat{n}(\xi h/2)$ is non-zero. To avoid the $n$th derivative of the product we make the following argument. We exclude the value $x=\xi h/2$ since then $r_{n}$ is equal to zero. Away from zero the derivative of the $\textrm{sgn}$ function is zero. Furthermore, it is $I_{1}$ that is difficult to represent accurately with polynomials and it contributes considerably more to the error than $F$ does, thus we keep only the dominant term obtained by the chain rule and write
\begin{align}
  r_n(x,y) \approx P_n(2y/h)\frac{(\alpha h)^{n}n!}{(2n)!}\textrm{sgn}\pars{\xi  h/2 - x}\alpha I_1^{(n)}\pars{\alpha(\xi  h/2 - x)}\frac{\pars{\xi  h/2 - x}\cdot\hat{n}(\xi h/2)}{|\xi  h/2 - x|}, \quad \xi\in[-1,1].
\label{eq:standard_interp_error_bound_dlp_2}
\end{align}

As for the single layer we set $\xi$ to zero, which results in the estimate
\begin{align}
  \Lambda_{D}(\alpha, n , \tilde{x}, h) &=
                  \frac{n!}{(2n)!}
                  (\alpha h)^{n+1}
                              I_1^{(n)}\pars{\frac{\alpha h }{2} |\tilde{x}|}\frac{-\mathrm{Im}(\tilde{x})}{|\tilde{x}|}\frac{1}{2}g_n(\tilde{x}),
\label{eq:lambdadefdbl}
\end{align}
where $g_n(\tilde{x})$ is given by \eqref{eq:leglogint} and  $\tilde{x} = 2x/h$ as before. Again we have reinserted the absolute value in the argument of the modified Bessel function, by removing the $\mathrm{sgn}$ function, since by \eqref{eq:defI} and  \eqref{eq:defIdiff} $I_{1}^{(n)}$ is an odd function for even $n$.

For the cancellation error, similarly
\begin{align}
  \Xi_{D}(\alpha,\tilde{x},h) = \epsilon_{\mathrm{mach}}\alpha h I_1\pars{\frac{\alpha h}{2} d}\frac{\mathrm{Im}(\tilde{x})}{|\tilde{x}|},
                    \quad d = 1 + |\tilde{x}|,
\label{eq:canc_err_est_3}
\end{align}
where $(y-x)\cdot \hat{n}(y)/|y-x|$ has been replaced with $\mathrm{Im}(\tilde{x})/|\tilde{x}|$.

\begin{figure}[htbp]
  \centering
      \includegraphics[width=0.49\textwidth]{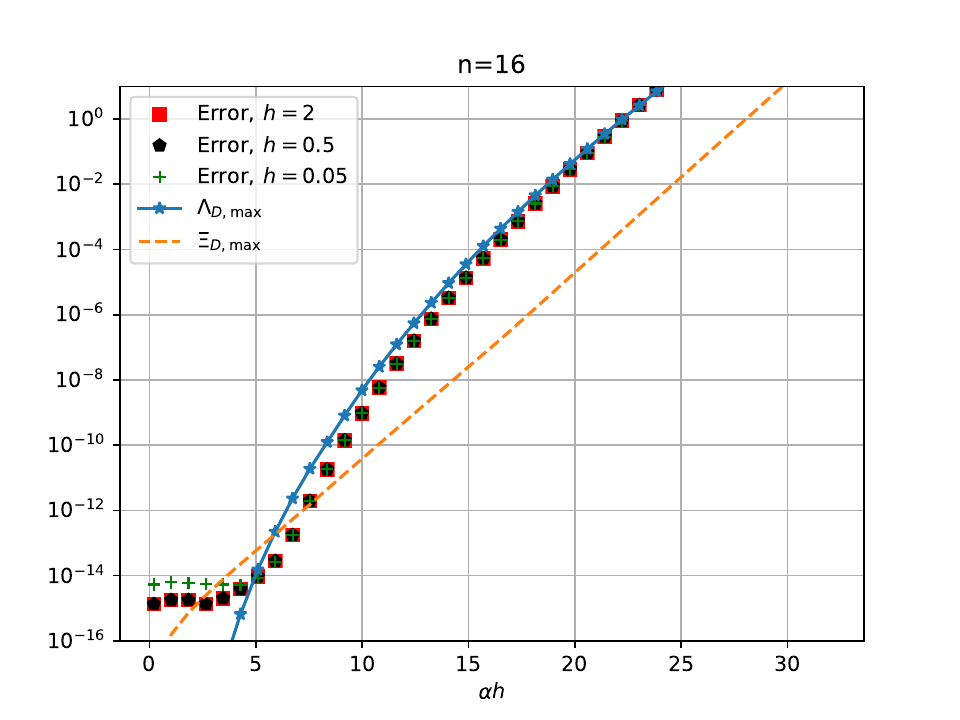}
      \\
  \includegraphics[width=0.49\textwidth]{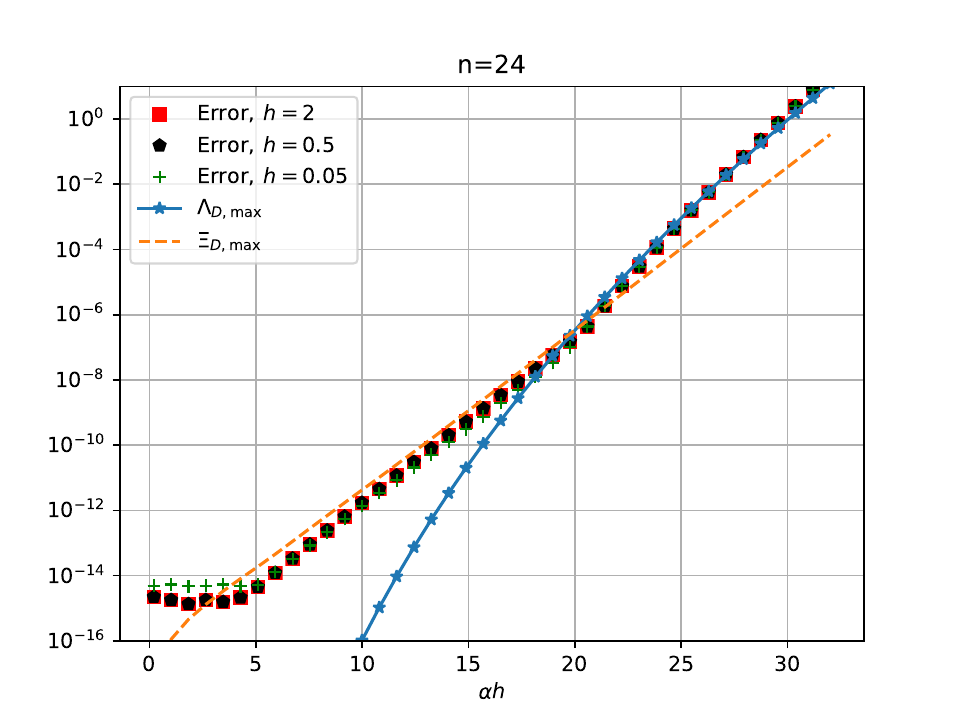}
  \hfill
    \includegraphics[width=0.49\textwidth]{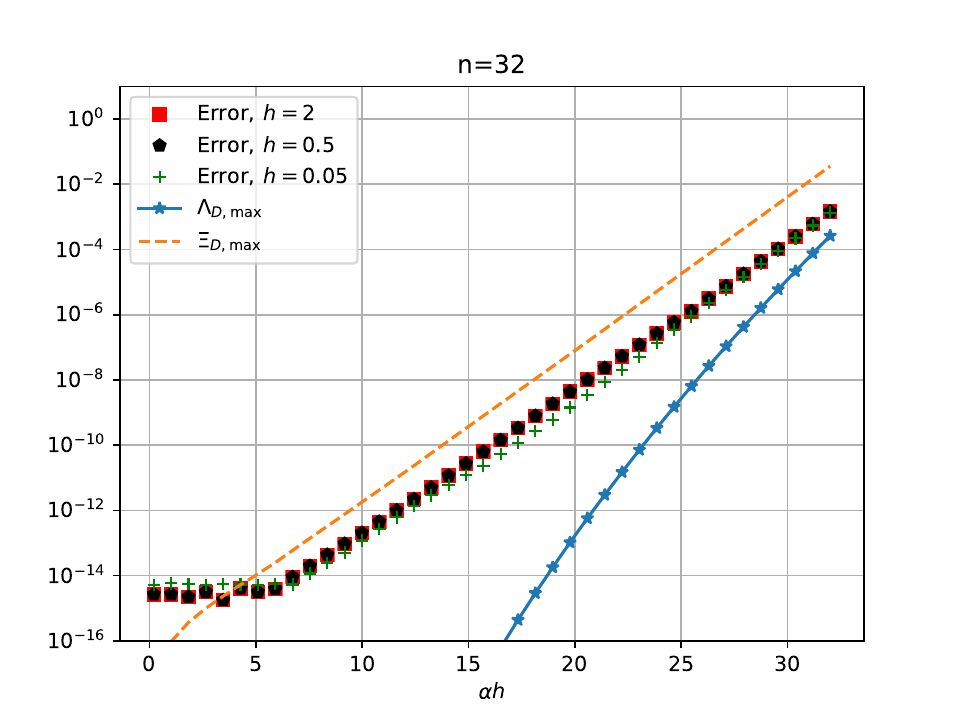}
  \caption{
Plots of measured actual total errors for different $h$, and error estimates $\Lambda_{D,\max}$ and $\Xi_{D,\max}$  for the target points in figure \ref{fig:errestn} for the double layer kernel for the modified Helmholtz equation \eqref{eq:modhelm_green_dbl}, with $n = 16$ (top), $n=24$ (left), and $32$ (right), over a range of $\alpha h$.}
  \label{fig:errestn2dbl}
\end{figure}

Clearly the error estimates depend purely on $\alpha h$, without any multiplying $h$-factor, in difference to the single layer kernel. For the double layer potential the error estimate is

\begin{equation}
  E_{D}(\alpha,n,h) = \tilde{E}_{D}(\alpha h, n) = \max \left(\Lambda_{D,\max},\Xi_{D,\max}\right),
\end{equation}
where $\Lambda_{D,\max}$ and $\Xi_{D,\max}$ are defined analogous to \eqref{eq:iperrmax} and \eqref{eq:ximax} for \eqref{eq:lambdadefdbl} and \eqref{eq:canc_err_est_3}. The predicted pure dependence on $\alpha h$ for the measured numerical error can be seen in figure \revMine{\ref{fig:errestn2dbl}}. The estimate for the cancellation error is not as tight as for the single layer kernel for $n=32$. Still, at its worst the  estimates are off by a digit of the actual error. Also, we see that it is possible to ignore the chain
rule and only differentiate $I_{1}$, but not $F$, and obtain good estimates. 

To derive error estimates for the kernels associated with the other modified PDEs a similar approach as above can be taken. The error estimate $E$ will have the form
\begin{equation}
\label{eq:errorscaled_general}
E(\alpha,n,h) = H(h)\tilde{E}(\alpha h,n),
\end{equation}
where simply $H(h)=h^p$ for some $p$.  

\begin{figure}[htbp]
  \centering
    \includegraphics[width=0.49\textwidth]{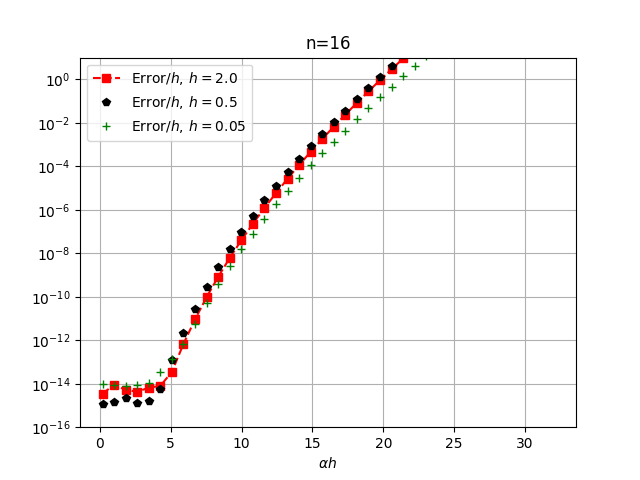}
  \caption{
Plot of measured actual total errors for different $h$, for the target points in figure \ref{fig:errestn} for the double layer kernel for the modified Stokes equation \eqref{eq:modstokes_green}, with $n=16$, over a range of $\alpha h$.}
  \label{fig:errestndblstokes}
\end{figure}

The discussion that we have put forward has given an understanding of the two sources of errors. Practically, only a value for $\Ceps$ needs to be determined; it sets panel lengths in the subdivision algorithm, introduced in the next section. Instead of deriving an error estimate, one can simply use figures such as figure \ref{fig:errestn}--\ref{fig:errestn2dbl} for the measured numerical errors for this purpose. Here, one needs to find the proper scaling $H(h)$ by which to scale the errors before plotting, such that they coincide for different $h$. In figure \revTwo{\ref{fig:errestndblstokes}} such results are presented for the double layer kernel for the modified Stokes equations \eqref{eq:modstokes_green}, with $p = 1$. The error curves do not collapse as neatly as for the modified Helmholtz equation, indicating that there is a scaling factor other than $h$. Still, the results can be used to set $\Ceps$, as demonstrated in section \ref{s:results}.

\subsection{Condition for $\alpha h$}
\label{ss:Cepsilon}
As was shown above, the error when evaluating the double layer
modified Helmholtz kernel with an $n$ point Gauss\revTwo{--}Legendre quadrature rule over
one panel with length $h$ only depends on $n$ and $\alpha h$. That
means that given the error tolerance $\varepsilon$, $\Ceps=\alpha h$
that yield this error can be read off from a plot such as in figure \ref{fig:errestn2dbl} for the given
$n$. Fixing that $\Ceps$, the condition (as written in \eqref{eq:alpha_h_crit}) is $\alpha
h \leq \Ceps$.

For other kernels, such as the single layer modified Helmholtz kernel,
there is not a pure $\alpha h$ dependence on the error, as we have
seen. However, to keep the procedure simple, we would suggest the following
approach. Pick a typical value of the panel length $h$ for the discretization at
hand. Then, divide the given tolerance $\varepsilon$ by this $h$, and
use this scaled tolerance when determining $\Ceps$ from plots of the
scaled error for the appropriate $n$, such as in figure \ref{fig:errestn} and figure \ref{fig:errestn2}.

In the subdivision algorithm, the reduction in error through this
extra factor of $h$ \revMine{of} refinement will then not be taken into account, and hence the
panels might be refined to be somewhat smaller than needed. But this
will allow us to keep the same simple structure of the subdivision algorithm for
all different kernels.

\section{Quadrature by recursive subdivision}
\label{s:algorithm}

To circumvent the problem of interpolatory quadrature failing for large $\alpha$, we here introduce an
algorithm for local refinement, based on panel subdivision. Given a single source panel
$\Gamma$, we assume that it is sufficiently short, relative to the
quadrature order $n$, for both the geometry and the layer density to be
well-represented by a polynomial, interpolated at the $n$ quadrature
points. We say that it is \emph{well-resolved}. For a given target
point $x$, we can then subdivide $\Gamma$ into a set of $M$ subpanels
$\{\Gamma_i\}_{i=1}^M$, interpolate our known quantities from the
quadrature nodes on $\Gamma$ to the quadrature nodes on those
subpanels, and then evaluate the layer potential at $x$ using the
subpanels. To ensure accuracy for a given tolerance, this subdivision is formed in a way
that guarantees that all subpanels either are short enough to satisfy
\eqref{eq:alpha_h_crit}, obtained via  \eqref{eq:estmax} or similarly, or are sufficiently far away from $x$,
relative to their own length, to not need kernel-split quadrature.

Before we can state our algorithm, a number of preliminaries are
needed.

\paragraph{Preimage of target}
Let $\gamma(t): \mathbb R \to \mathbb C$ be the mapping from the
standard interval $[-1,1]$ to the panel $\Gamma$. Then $z$, such that
$\gamma(z)=x$, is the \emph{preimage} of the target point $x$. The
preimage $z$ is real-valued if $x\in\bdry$, and complex-valued
otherwise. We here assume that we know the value of $z$, but $\gamma(t)$ need not be a known function; see
\cite{AfKlinteberg2018} for a discussion on how construct a numerical representation.

\paragraph{Subpanels and subintervals}
A subdivision of $\Gamma$ is defined by a set of edges in the
parametrization, $\{-1=t_1, t_2, \dots, t_{M+1}=1\}$, such that a
subpanel $\Gamma_i$ is given by the mapping of the subinterval
$[t_i, t_{i+1}]$ under $\gamma$. We can, by a linear scaling, define
the \emph{local mapping} that maps the standard interval $[-1,1]$ to
$\Gamma_i$ as
\begin{align}
  \gamma_i(t)=\gamma\pars{t_i +  \frac{\Delta t_i}{2}(t+1)},
\end{align}
where $\Delta t_i=t_{i+1}-t_i$. Given the preimage $z$, the
\emph{local preimage} $z_i$, such that $\gamma_i(z_i)=x$, is given by
\begin{align}
  z_i = \frac{2}{\Delta t_i}(z-t_i) - 1 .
  \label{eq:local_preimage}
\end{align}

\paragraph{Near evaluation criterion}
Given the preimage $z$ of a point $x$ close to a panel (or subpanel)
$\Gamma$, it is possible to compute an accurate estimate of the
quadrature error when evaluating the layer potential using $n$-point
Gauss\revTwo{--}Legendre quadrature. Detailed discussions can be found in
\cite{AfKlinteberg2016quad,AfKlinteberg2018}. To leading order, the
error is proportional $\rho(z)^{-2n}$, where $\rho$ is the elliptical
radius of the Bernstein ellipse on which $z$ lies,
\begin{align}
  \rho(z) = \abs{z + \sqrt{z^2-1}},
  \label{eq:bernstein_radius}
\end{align}
where $\sqrt{z^2-1}$ is defined as $\sqrt{z+1}\sqrt{z-1}$ with $-\pi < \arg\pars{z\pm 1}\leq \pi$. For a given
kernel $G$ and error tolerance $\varepsilon$, it is then possible to
introduce a cutoff radius $\Reps$, such that kernel-split quadrature
must be used for
\begin{align}
  \rho(z) < \Reps,
  \label{eq:near_eval_crit}
\end{align}
and otherwise Gauss\revTwo{--}Legendre quadrature is sufficiently accurate. We refer to this as the target point and the source points being \emph{well-separated}, otherwise they are considered to be \emph{close}. See \cite{AfKlinteberg2018} on how to set $\Reps$ for a given tolerance;  we use $\Reps= 3.5$ for $n = 16$, which corresponds to a tolerance of $10^{-14}$.

Later we will use that the inverse of $\rho$ has a particularly simple form  in the
special case when $z$ lies on the imaginary axis,
\begin{align}
  \begin{split}
    z &= \pm \im b, \quad b>0,\\
    \rho(z) &= b+\sqrt{b^2+1}, \\
    z(\rho) &= \pm \frac{\im (\rho^2 - 1)}{2\rho} .
  \end{split}
        \label{eq:bernstein_center}
\end{align}

\paragraph{Interpolation and upsampling}

To interpolate data from the $n$ original Gauss\revTwo{--}Legendre nodes on
$[-1,1]$, to $m$ new Gauss\revTwo{--}Legendre nodes on a subinterval
$[t_i, t_{i+1}] \subset [-1,1]$, we use barycentric Lagrange interpolation
\cite{Berrut2004}. By \emph{upsampling} we refer to the special case
of interpolating from $n$ to $2n$ Gauss\revTwo{--}Legedre nodes, both on
$[-1,1]$.

\paragraph{Subinterval length criterion}

When a new subpanel $\Gamma_i$ is formed on $\Gamma$, we need to check
if it satisfies the kernel-split accuracy criterion
\eqref{eq:alpha_h_crit}, which requires knowledge of the arc length of
the subpanel, denoted $h_i$. Assuming that $\gamma'(t)$ does not vary
rapidly on $\Gamma$, a good approximation to $h_i$ is
$h_i \approx h\Delta t_i/2$, where $h$ is the arc length of
$\Gamma$. We can now combine this approximation with
\eqref{eq:alpha_h_crit}, to get an accuracy criterion formulated in
subinterval size,
\begin{align}
  \Delta t_i \le \frac{2 \Ceps}{\alpha h}.
  \label{eq:alpha_dt_crit}
\end{align}
In particular, this allows us to write down the maximum length of
subintervals on which product integration can be used,
\newcommand{\dtmax}{{\Delta t_{\max}}}
\begin{align}
  \Delta t_{\max} = \frac{2 \Ceps}{\alpha h} .
  \label{eq:dt_max}
\end{align}
Here, $\Ceps$ can be obtained via \eqref{eq:estmax}, as a precomputation step.

\subsection{Algorithm}

Our algorithm proceeds with creating a division of $[-1,1]$ into
subintervals, which corresponds to a division of $\Gamma$ into
subpanels.

For a target point $x$ with preimage $z$ such that $\Re z \in (-1,1)$,
the first step is to create a subinterval centered on $\Re z$, with
length set to satisfy both of the conditions \eqref{eq:near_eval_crit}
and \eqref{eq:dt_max}. The centering ensures that the subpanels will
not introduce new edges or quadrature nodes that are close enough to
$z$ to degrade precision.

If this initial subinterval has length $\Delta t_c$, then the preimage
of $x$ in that local frame will be $z_c = \im b_c$, with
$b_c=2 \Im z / \Delta t_c$ and \eqref{eq:bernstein_center} is applicable.

If we wish the local preimage $z_c$ to be
just beyond the limit where kernel-split quadrature is needed, then we
must set $\Delta t_c$ such that $\rho(z_c) = \Reps$. From
\eqref{eq:bernstein_center}, we can derive that this is satisfied when
$\Delta t_c = \Delta t_{\text{direct}}$,
\begin{align}
  \Delta t_{\text{direct}} = |\Im z|  \frac{4\Reps}{\Reps^2 - 1} .
  \label{eq:dt_direct}
\end{align}
For the subinterval to be contained within $[-1,1]$, it may not be
bigger than twice the distance between $\Re z$ and the closest interval edge,
\begin{align}
  \Delta t_{\text{edge}} = 2(1-|\Re z|).
\end{align}
Now we set the initial subpanel as large as possible, while still
ensuring that the quadrature from it is accurate, and that it falls
within $[-1,1]$,
\begin{align}
  \Delta t_c = \min\pars{\Delta t_{\text{edge}}, \max\pars{\Delta t_{\text{direct}}, \dtmax} }.
  \label{eq:deltatc}
\end{align}
Here, $\Delta t_{\text{direct}}$ depends on the location of the target point, while $\dtmax$ does not. Thus, for points sufficiently far away $\dtmax$ might be larger than $\Delta t_{\text{direct}}$.

By \eqref{eq:deltatc} we have the initial subdivision
$\{-1, \Re z-\Delta t_c/2, \Re z+\Delta t_c/2, 1\}$. The center
subinterval is now acceptable, and we proceed by recursively bisecting
each remaining subinterval until either its length
satisfies \eqref{eq:alpha_dt_crit}, or the local preimage of $x$
satisfies \eqref{eq:near_eval_crit}.

For target points such that $\Re z \notin (-1,1)$, we can skip the
process of carefully selecting the length of the nearest subinterval,
and proceed immediately with recursive bisection of $\{-1,1\}$. This
completes the algorithm, which we list in its entirety in algorithm
\ref{alg:subdivide}.

\begin{algorithm}[htbp]
  \caption{Given a panel of length $h$ and a nearby target point with
    preimage $z$, create a subdivision of $[-1,1]$ that allows the
    layer potential to be accurately evaluated, using either direct
    Gauss\revTwo{--}Legendre quadrature or kernel-split quadrature.}
  \label{alg:subdivide}
  \begin{algorithmic}
    \Function{create\_subdivision}{$z, h, \alpha, \Ceps, \Reps$}
    \State \revTwo{$\Delta t_{\max} \gets 2 \Ceps / (\alpha h)$}
    \If{$\abs{\Re z} \ge 1$} \Comment Preimage outside interval, recursively bisect all of it.
    \State \Return \Call{recursive\_bisection}{$-1, 1,\Reps, \dtmax, z$}.
    \Else
    \State $\Delta t_{\text{direct}} \gets 4 |\Im z| \Reps / (\Reps^2 - 1)$
    \State $\Delta t_{\text{edge}} = 2(1-|\Re z|)$
    \State $\Delta t_c \gets  \min\pars{\Delta t_{\text{edge}}, \max\pars{\Delta t_{\text{direct}}, \dtmax} }$
    \State $t_a \gets \Re z - \Delta t_c/2$
    \State $t_b \gets \Re z + \Delta t_c/2$
    \State \Comment Center interval is now acceptable, recursively bisect remainder subintervals.
    \State $S_1 \gets $ \Call{recursive\_bisection}{$-1, t_a,\Reps, \dtmax, z$}
    \State $S_2 \gets  \{t_a, t_b\}$
    \State $S_3 \gets $ \Call{recursive\_bisection}{$t_b, 1,\Reps, \dtmax, z$}
    \State \Return $S_1 \cup S_2 \cup S_3$
    \EndIf
    \EndFunction
    \\
    \Function{recursive\_bisection}{$t_1, t_2,\Reps, \dtmax, z$}
    \If{$t_1 < t_2$}
    \State $\Delta t_{\text{sub}} \gets t_2-t_1$
    \State $z_{\text{sub}} \gets 2(z-t_1)/\Delta t_{\text{sub}} - 1$   \Comment From \eqref{eq:local_preimage}.
    \If{$\rho(z_{\text{sub}}) < \Reps$ \textbf{and} $\Delta t_{\text{sub}} > \dtmax$}
    \Comment Using \eqref{eq:bernstein_radius}.
    \State \Comment Kernel-split must be used, but interval still too large. Continue bisection.
    \State $t_{\text{mid}} \gets t_1 + \Delta t_{\text{sub}}/2$
    \State $S_1 \gets $ \Call{recursive\_bisection}{$t_1, t_{\text{mid}},\Reps, \dtmax, z$}
    \State $S_2 \gets $ \Call{recursive\_bisection}{$t_{\text{mid}}, t_2,\Reps, \dtmax, z$}
    \State \Return $S_1 \cup S_2$
    \EndIf
    \EndIf
    \State \Return $\{t_1, t_2\}$ \Comment Subinterval passed.
    \EndFunction
  \end{algorithmic}
\end{algorithm}

\section{Numerical results}
\label{s:results}
To test the robustness of our method across a range of $\alpha$
values, we solve the modified Helmholtz equation and the modified Stokes equations. Based on the discussion in section \ref{s:errorintro} the parameter $\Ceps$ is set to confirm that prescribed error tolerances are satisfied.
\subsection{The modified Helmholtz equation}
\begin{figure}[htb]
  \centering
  \subfloat[$K_{0}(10\|x-x_{0}\|)$]{
    \includegraphics[width=0.45\textwidth,trim={4cm 8.5cm 6cm 5cm},clip]{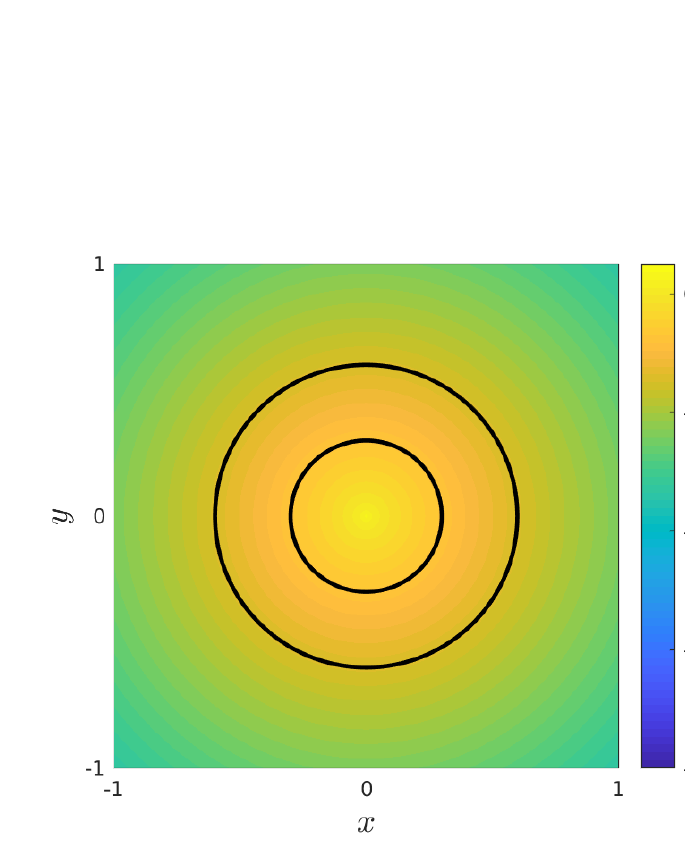}
  }
  \hfill
  \subfloat[$K_{0}(100\|x-x_{0}\|)$]{
    \includegraphics[width=0.45\textwidth,trim={4cm 8.5cm 6cm 5cm},clip]{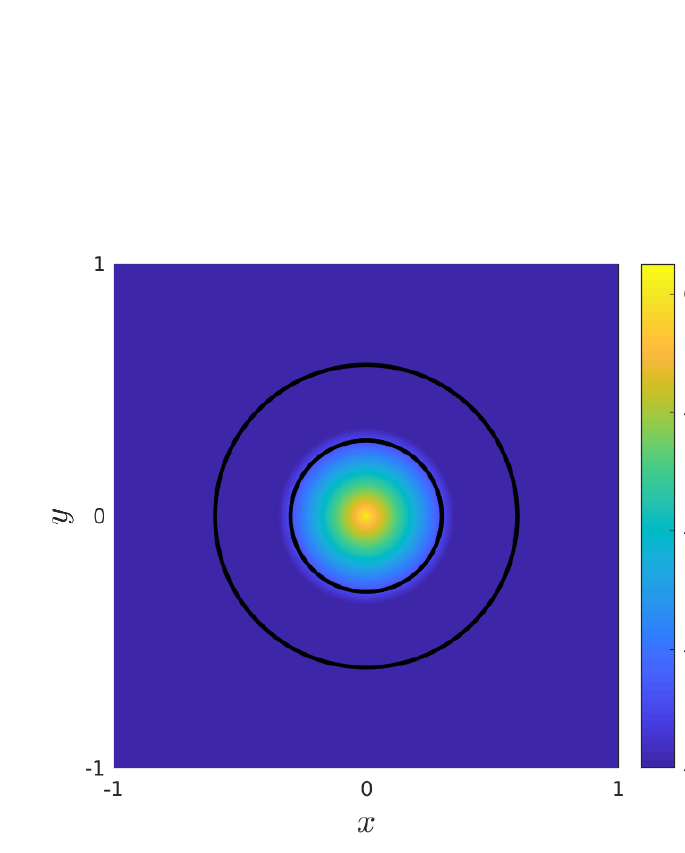}
  }
  \caption{Plots of $K_{0}(\alpha\|x - x_{0}\|)$, the modified Bessel function of the second kind of zeroth order, for $\alpha = 10,100$ and $x_0 =  (0.01, 0.01)^T$. The black circles of radius $0.3$ and $0.6$ define the annulus that is computational domain.}
  \label{fig:modhelmsolution}
\end{figure}

We solve the modified Helmholtz equation \eqref{eq:modhelm} inside the annulus centered at the origin, defined by a
circle of radius 0.3, and a circle of radius 0.6, with Dirichlet
boundary conditions given by a fundamental solution
\eqref{eq:modhelm_green} centered in the inner circle,
\begin{align}
  (\Delta - \alpha^2)u &= 0,   &x &\in \domain, \\
  \label{eq:solmodhelm}
  u &= K_0(\alpha|x-x_0|), \quad  &x &\in \bdry, \\
  x_0 &= (0.01, 0.01)^T. &&
\end{align}
The exact solution to this problem is equal to the expression for the Dirichlet boundary
condition, evaluated in $\domain$. See figure \ref{fig:modhelmsolution} for a visualization of \eqref{eq:solmodhelm} for $\alpha = 10,\,100$. We see that $K_0(\alpha|x-x_0|)$ decays rapidly for $\alpha = 100$.

To solve the above problem numerically, we represent the
solution using the double layer potential
\begin{align}
  u(x) = \int_\bdry G(x,y)\sigma(y)  \dif S(y),
  \label{eq:modhelm_dbl}
\end{align}
where $G$ is the double layer kernel \eqref{eq:modhelm_green_dbl}.

Enforcing the boundary condition \revMine{\eqref{eq:solmodhelm}} gives a second kind integral equation in $\sigma$,
\begin{align}
  \sigma(x) + \int_\bdry  G(x, y) \sigma(y)  \dif S(y)
  = K_0(\alpha|x-x_0|), \quad x \in \bdry.
  \label{eq:biemodhelmslp}
\end{align}
We solve this using the Nystr\"om method, discretizing the boundary
using $16$-point Gauss\revTwo{--}Legendre panels, with $15$ panels on the inner
circle, and $30$ panels on the outer. For the bound
\eqref{eq:dt_direct} we use use $\Reps = 3.5$, and for \eqref{eq:alpha_dt_crit} we set $\Ceps = 3.7$ by reading off figure \ref{fig:errestn2dbl} to achieve a tolerance of $10^{-14}$. Following the
notation of \eqref{eq:split}, the kernel-split of
\eqref{eq:modhelm_green_dbl} is given by  \eqref{eq:modhelm_dbl_G0} -- \eqref{eq:modhelm_dbl_GC_domain}.
\begin{figure}[htb]
  \centering
  \subfloat[Error when evaluating the layer potential close to the
boundary, relative error in max-norm computed as $\|u -
u_{\text{ref}}\|_\infty / \|\sigma\|_\infty$.]{
    \includegraphics[width=0.45\textwidth]{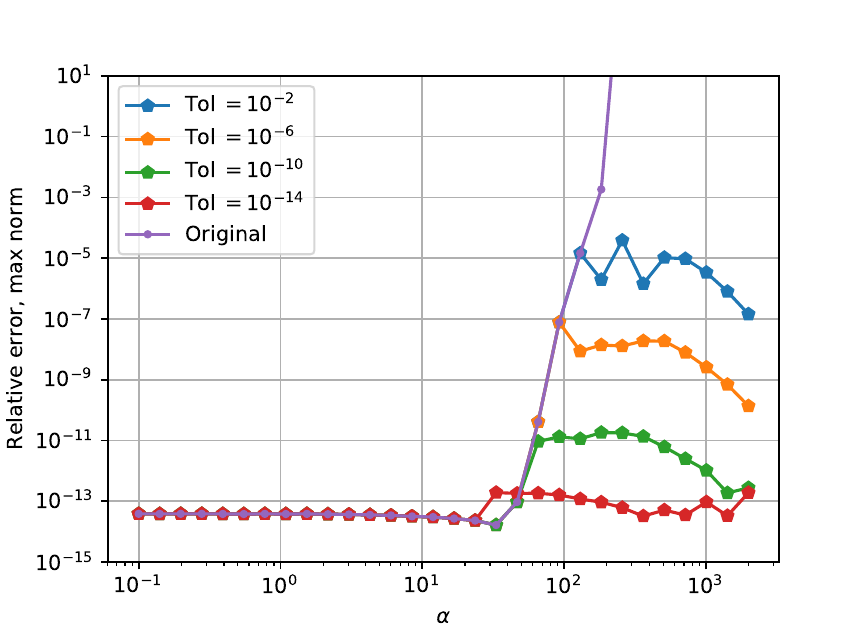}
  }
  \hfill
  \subfloat[Time to assemble the matrix blocks that correspond to
neighboring panels, which are where the $\log$-correction need to be
added for the tolerance $10^{-14}$.]{
    \includegraphics[width=0.45\textwidth]{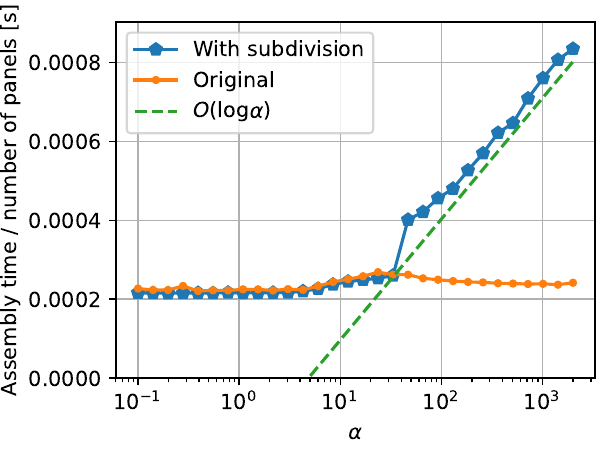}
  }
  \caption{Comparison of the original kernel-split algorithm, denoted
    ``Original'', and our adaptive algorithm, denoted ``Tol $=$ tolerance
    subdivision'' for given tolerances, when solving our test problem for the modified Helmholtz equation for a large range of
    $\alpha$. We test the solution up to $\alpha \approx 2000$; for
    larger values of $\alpha$ the solution is about round-off in the
    entire domain.}
  \label{fig:vary_alpha}
\end{figure}

For a large range of $\alpha$ values, we solve the integral equation,
and then evaluate the solution at 15 random points on a circle of
radius 0.301 (very close to the inner boundary). Quadrature by recursive bisection is applied both when solving for the density \eqref{eq:biemodhelmslp} and for evaluating the solution \eqref{eq:modhelm_dbl}.

 The results, shown in figure
\ref{fig:vary_alpha}, demonstrate that our subdivision algorithm is
capable of avoiding the catastrophic loss of accuracy otherwise
present above a threshold $\alpha$, and that the additional cost
incurred from it is proportional to $\log\alpha$. We do not satisfy the error tolerance of $10^{-14}$; around one digit of
accuracy appears to be lost, presumably due to the additional interpolation steps
involved.

We also try other values for $\Ceps$ to confirm that we stay under a given tolerance. Based on figure \ref{fig:errestn2dbl} we have $\Ceps = 18.7,\,12.8,\,8.2$ for the tolerances $10^{-2},\,10^{-6},\,10^{-10}$. The largest corresponding relative error for each tolerance, taken over the plotted range of values of $\alpha$, are $3.7\cdot 10^{-5},\,6.7\cdot 10^{-8},\,1.8\cdot 10^{-11}$. Clearly the errors are at least one digit below the prescribed tolerance. It is not suprising that the set value for $\Ceps$ gives a lower error than figure \ref{fig:errestn2dbl} suggests. In the subdivision algorithm introduced in this paper, the first panel in the recursive scheme always has as its center the projection of the target point on the boundary. Therefore target points are never close to panel edges, where errors from the kernel-split quadrature tend to be greater than for target points towards the panel's center \cite{Helsing2008}.

\subsection{The modified Stokes equations}
We solve the modified Stokes equations \eqref{eq:modstokes} in a setting analogous to \eqref{eq:solmodhelm}; the geometry, and the discretization thereof, is the same, and the boundary data and the solution is given by the associated fundamental solution \eqref{eq:modstokes_green}. As opposed to the modified Helmholtz equation, no error estimates are presented in this paper for the modified Stokes equation, which can be used to set $\Ceps$. Instead, one can use the measured numerical errors for a flat panel for different values of $h$, each scaled with $1/h^p$, as suggested in section \ref{ss:otherkernels}. Here $p$ is chosen such that the error curves collapse. For the modified Stokes in a double layer formulation, such results are shown in figure \ref{fig:errestndblstokes}, with $p = 1$.

\begin{figure}[htbp]
  \centering
    \includegraphics[width=0.49\textwidth]{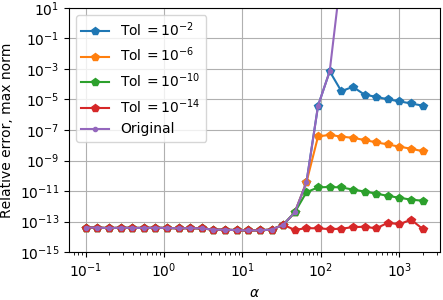}
  \caption{
Comparison of the original kernel-split algorithm, denoted ``Original'', and our adaptive algorithm, denoted ``Tol = tolerance subdivision'' for given tolerances, when solving our test problem for the modified Stokes equations for a large range of $\alpha$.  
}
  \label{fig:vary_alpha_stokes}
\end{figure}

Following the methodology presented in section \ref{ss:Cepsilon}, we set $\Ceps$ to achieve certain tolerances. For the tolerances $10^{-2},\,10^{-6},\,10^{-10},\,10^{-14}$, we set $\Ceps = 16.7,\,11.5,\,7.4,\,4.4$ based on figure \ref{fig:errestndblstokes}. The corresponding \revMine{maximum errors over all alpha} are $3.8\cdot 10^{-5}, \,4.8\cdot 10^{-8},\,1.1\cdot 10^{-11},\,5.3\cdot 10^{-14}$, \revMine{shown} in figure \ref{fig:vary_alpha_stokes}. The magnitudes of the errors, for a given tolerance, is consistent with the results for the modified Helmholtz equation, presented in the previous section. The errors are well below the prescribed tolerance, meaning more subdivisions are applied than necessary. The hypothesis is the same as for the modified Helmholtz; for the subdivision algorithm the target point, or the projection thereof, on the boundary is always centered on the new panel. Thus its never close to the panel edges, where the error tend\revMine{s} to be greater.

\section{Conclusions}
\label{s:conclusions}

We present a robust recursive algorithm that allows the method of
Helsing et al. to be applied, for any $\alpha$, to the modified
Helmholtz equation, modified biharmonic equation and modified Stokes
equations \revOne{on smooth geometries}.  Before, this was not possible for large $\alpha$, which
corresponds to small time-steps with semi-implicit marching schemes for
the heat equation and the time-dependent Stokes and Navier-Stokes
equations.  Our algorithm is fully adaptive, and the additional
computational time it requires scales as $\log \alpha$. Our choice of
the parameters $\Ceps$ and $\Reps$ is based on numerical observations
and provide excellent results. 
\section*{Acknowledgments}

The authors gratefully acknowledge the support from
the Knut and Alice Wallenberg Foundation under grant no. 2016.0410
(L.a.K.),
and by the Swedish Research Council under Grant No. 2015-04998 (F.F. and A.K.T.)
\section*{Declarations}
The authors declare no competing interests.
\clearpage
\appendix
\begin{appendices}

\section{Formulas for recursive computation of the Legendre-log
    integrals}
\label{sec:leglog}

\subsection{Result}
The integrals
\begin{align}
  g_n(x) = \int_{-1}^1 P_n(\genvar) \log(\genvar-x) \dif \genvar, \quad n=0,1,2,\dots ,
  \label{eq:leglog_def}
\end{align}
can be computed recursively using the formulas
\begin{align}
  \begin{split}
    g_0(x) &=
    (1 - x)\log(1 - x) +
    (1 + x) \log(-1-x) - 2,
    \\
    g_1(x) &=
    \frac{1}{2} \left(
      \left(1-x^2\right) \log (1-x)
      -
      \left(1-x^2\right) \log (-x-1)
      -2 x
    \right),
    \\
    g_2(x) &= \frac{1}{3}\pars{3xg_1 + 2}, \\
    &\vdots \\
    g_{n+1}(x) &= \frac{1}{n+2}\pars{(2n+1)xg_{n}(x) - (n-1)g_{n-1}(x)},
    \quad n\ge2.
  \end{split}
              \label{eq:leglog_recur}
\end{align}
The first two formulas, for $g_0(x)$ and $g_1(x)$, have finite limits
as $x \to \pm 1$.

\subsection{Proof}
The formulas for $g_0$ and $g_1$ follow from direct integration of
\eqref{eq:leglog_def}, with $P_0(\genvar)=1$ and $P_1(\genvar)=\genvar$. For $n>1$, the
following two results are useful:
\begin{align}
  (n+1) P_{n+1}(\genvar) &= (2n+1) \genvar P_n(\genvar) - n P_{n-1}(\genvar),
                     \label{eq:leg_recur}
  \\
  (2n+1) P_n(\genvar) &= \dod{}{\genvar} \pars{ P_{n+1}(\genvar) - P_{n-1}(\genvar) } .
                  \label{eq:leg_diffident}
\end{align}
Insertion of the recursion formula \eqref{eq:leg_recur} into
\eqref{eq:leglog_def} gives
\begin{align}
  (n+1)g_{n+1}(x)
  &= (2n+1) \int_{-1}^1
    \genvar P_n(\genvar)
    \log(\genvar-x) \dif \genvar
    -
    n \int_{-1}^1
    P_{n-1}(\genvar)
    \log(\genvar-x) \dif \genvar \\
  &= (2n+1) \int_{-1}^1
    P_n(\genvar)(\genvar-x+x)
    \log(\genvar-x) \dif \genvar
    -
    n g_{n-1}(x)
  \\ &= (2n+1) \int_{-1}^1
    P_n(\genvar)(\genvar-x)
       \log(\genvar-x) \dif \genvar
       +
       (2n+1) x g_n(x)
       -
       n g_{n-1}(x) .
\end{align}
Integration by parts of remaining integral, using \eqref{eq:leg_diffident} and $((\genvar-x)\log(\genvar-x))'=1+\log(\genvar-x)$,
\begin{align}
  &(2n+1)\int_{-1}^1
  P_n(\genvar)(\genvar-x)
  \log(\genvar-x) \dif \genvar  =
    \\
  &=\underbrace{
    \left[
  \pars{P_{n+1}(\genvar) - P_{n-1}(\genvar)}\log(\genvar-x)
    \right]_{-1}^1
    }_{I_1}
  -
    \underbrace{
    \int_{-1}^1
  \pars{P_{n+1}(\genvar) - P_{n-1}(\genvar)}\pars{1+\log(\genvar-x)}
  \dif \genvar
    }_{I_2} .
\end{align}
We have $P_n(1)=1$ and $P_n(-1)=(-1)^n$, so it follows that $I_1=0$. Furthermore,
\begin{align}
  \int_{-1}^1 P_n(\genvar) \dif \genvar = 2\delta_{0n} =
  \begin{cases}
    2, & \text{if} \quad n=0\\
    0, & \text{if} \quad n\ge 1 .
  \end{cases}
\end{align}
Consequently, since $n \ge 0$,
\begin{align}
  I_2 &= -2\delta_{0(n-1)} + g_{n+1}(x) - g_{n-1}(x),
  \\
  \Rightarrow (n+1)g_{n+1}(x)
  &=
    2\delta_{0(n-1)} - g_{n+1}(x) + g_{n-1}(x)
    +
    (2n+1) x g_n(x)
    -
    n g_{n-1}(x),
  \\
\Rightarrow  (n+2)g_{n+1}(x)
  &=
    (2n+1) x g_n(x)
    -
    (n-1) g_{n-1}(x)
    +
    2\delta_{0(n-1)},
\end{align}
which leads to the recursion formulas \eqref{eq:leglog_recur}.

\section{Kernel-splits}
\label{sec:kernelsplits}
Here we complement the kernel-split for the modified Helmholtz equation with kernel-splits for the modified biharmonic equation and the modified Stokes equations.

\subsection{Modified biharmonic equation}
The single layer kernel for the modified biharmonic equation is
\begin{equation}
  G(x,y) = \frac{-1}{2\pi\alpha^2}(\log\revTwo{|y-x|} + K_{0}(\alpha \revTwo{|y-x|})).
\end{equation}
Inserting \eqref{eq:K0_split} the explicit split becomes
\begin{align}
  G^{S} &= \frac{-1}{2\pi\alpha^{2}} ( K_{0}^{S}(\alpha\revTwo{|y-x|}) + I_{0}(\alpha\revTwo{|y-x|})\log\alpha),\\
  \label{eq:mbslpgl}
  G^{L} &= \frac{-1}{2\pi\alpha^{2}} ( 1 + I_{0}(\alpha\revTwo{|y-x|})),\\
  G^{C} &= 0.
\end{align}
Here, $G^{L}$ contains the modified Bessel function $I_{0}$, like the single layer kernel for the modified Helmholtz equation. The double layer kernel is
\begin{align}
  G(x,y) &= \frac{\partial}{\partial\hat{n}(y)}\frac{-1}{2\pi\alpha^2}(\log\revTwo{|y-x|} + K_{0}(\alpha \revTwo{|y-x|}))
           = \frac{-1}{2\pi\alpha^2}\left(\frac{1}{\revTwo{|y-x|}} + \alpha K_{1}(\alpha \revTwo{|y-x|})\right) \frac{(y-x)\cdot \hat{n}(y)}{\revTwo{|y-x|}}\\
  &=  \frac{-1}{2\pi\alpha^2}\left(\frac{1}{\revTwo{|y-x|}} + \alpha K_{1}^{S}(\alpha \revTwo{|y-x|}) + \frac{1}{\revTwo{|y-x|}} + \alpha I_{1}(\alpha \revTwo{|y-x|})\log(\alpha \revTwo{|y-x|})\right) \frac{(y-x)\cdot \hat{n}(y)}{\revTwo{|y-x|}}
\end{align}
by \eqref{eq:K1_split}. The resulting decomposition is
\begin{align}
  G^{S} &= \frac{-1}{2\pi\alpha} ( K_{1}^{S}(\alpha\revTwo{|y-x|}) + I_{1}(\alpha\revTwo{|y-x|})\log\alpha) \frac{(y-x)\cdot \hat{n}(y)}{\revTwo{|y-x|}},\\
  G^{L} &= \frac{-1}{2\pi\alpha}I_{1}(\alpha\revTwo{|y-x|}) \frac{(y-x)\cdot \hat{n}(y)}{\revTwo{|y-x|}},\\
  G^{C} &= \frac{-1}{\pi\alpha^2}.
\end{align}

\subsection{The modified Stokes equations}
We present only the double layer kernel for the modified Stokes equations, also known as the stresslet. In the Einstein summation convention it has the closed-form
\begin{equation}
  \label{eq:modstokes_green}
  G_{ijk} = \alpha^2 \mathcal{G}_{1}(\alpha \|\mathbf{r}\|)(\delta_{jk}r_{i} + \delta_{ik}r_{j} + \delta_{ij}r_{k}) + \alpha^{4}\mathcal{G}_{2}(\alpha \|\mathbf{r}\|)r_{i}r_{j}r_{k} + \alpha^{2}\mathcal{G}_{3}(\alpha \|\mathbf{r}\|)\delta_{ik}r_{j},
\end{equation}
where $\delta_{ij}$ is the Kronecker delta and $\mathbf{r} = \mathbf{x}-\mathbf{y}$. Here, the functions $\mathcal{G}_{1}$--$\mathcal{G}_{3}$ are
\begin{align}
  \mathcal{G}_{1}(\genvar) &= -\frac{ 2\genvar^2 K_{0}(\genvar) +
                             (\genvar^2 + 4)\genvar K_{1}(\genvar)-4}{ 2\pi \genvar^4}, \\
  \mathcal{G}_{2}(\genvar) &= \frac{ 4\genvar^2 K_{0}(\genvar) +
                             (\genvar^2 + 8)\genvar K_{1}(\genvar)-8}{ \pi \genvar^6}, \\
  \mathcal{G}_{3}(\genvar) &= \frac{\genvar K_{1}(\genvar)-1}{ 2 \pi \genvar^2}.
\end{align}

Since they are expressed in terms of the modified Bessel functions $K_{0}$ and $K_{1}$ we can use the decompositions \eqref{eq:K0_split} and  \eqref{eq:K1_split} to write the expressions in explicit singularities. We have
  
\begin{align}
  \mathcal{G}_{1}(\genvar) &=  \mathcal{G}_{1}^{S}(\genvar) +  \mathcal{G}_{1}^{L}(\genvar)\log \genvar,\\
  \mathcal{G}_{2}(\genvar) &=  \mathcal{G}_{2}^{S}(\genvar) +  \mathcal{G}_{2}^{L}(\genvar)\log \genvar + \frac{1}{8\pi \genvar^{2}}-\frac{1}{\pi \genvar^4},\\
  \mathcal{G}_{3}(\genvar) &=  \mathcal{G}^{3}_{S}(\genvar) +  \mathcal{G}_{3}^{L}(\genvar)\log \genvar,
\end{align}
where
\begin{align}
  \mathcal{G}_{1}^{S} &= -\frac{2\genvar K_{0}^{S}(\genvar)+(\genvar^2
                        + 4)K_{1}^{S}(\genvar) + \genvar}{2\pi
                        \genvar^{3}},\quad &\mathcal{G}_{1}^{L} &=
                                                                  \frac{2\genvar
                                                                  I_{0}(\genvar) - (\genvar^2+4)I_{1}(\genvar)}{2\pi \genvar^{3}},\\
  \mathcal{G}_{2}^{S} &= \frac{32\genvar
                        K_{0}^{S}(\genvar)+8(\genvar^2 + 8)
                        K_{1}^{S}(\genvar) - \genvar(\genvar^2 -
                        16)}{8\pi \genvar^{5}},\quad
                                           &\mathcal{G}_{2}^{L} &=
                                                                  \frac{(\genvar^{2}
                                                                  +
                                                                  8)I_{1}(\genvar)
                                                                  -
                                                                  4\genvar
                                                                  I_{0}(\genvar)}{\pi \genvar^{5}},\\
   \mathcal{G}_{3}^{S} &= \frac{K_{1}^{S}(\genvar)}{2\pi \genvar},\quad & \mathcal{G}^{3}_{L} &=\frac{I_{1}(\genvar)}{2\pi \genvar}.\\
\end{align}
With these definitions the kernel-split form of the stresslet is
\begin{equation}
 T_{ijk}(\mathbf{r}) = \left( T_{ijk}^{S}(\mathbf{r}) +  T_{ijk}^{L}(\mathbf{r})\log(\alpha)\right) + T_{ijk}^{L}(\mathbf{r})\log\|\mathbf{r}\|  +  T_{ij}^{C}(\mathbf{r})\frac{r_{k}}{\|\mathbf{r}\|} + T^{Q}\frac{r_{i}r_{j}r_{k}}{\|\mathbf{r}\|^{4}},
\end{equation}
\begin{equation}
  G_{ijk}(\mathbf{r}) = G_{ijk}^{S}(\mathbf{r}) +  T_{ijk}^{L}(\mathbf{r})\log\|\mathbf{r}\| +  T_{ij}^{C}(\mathbf{r})\frac{r_{k}}{\|\mathbf{r}\|} + T^{Q}\frac{r_{i}r_{j}r_{k}}{\|\mathbf{r}\|^{4}}
\end{equation}
where
\begin{align}
  G^{S}_{ijk}(\mathbf{r}) &= \alpha^{2}\mathcal{T}_{1}^{S}(\alpha \|\mathbf{r}\|)(\delta_{jk}r_{i}+\delta_{ik}r_{j}+\delta_{ij}r_{k}) + \alpha^{4}\mathcal{T}_{2}^{S}(\alpha \|\mathbf{r}\|)r_{i}r_{j}r_{k}  + \alpha^{2}\mathcal{T}_{3}^{S}(\alpha \|\mathbf{r}\|)\delta_{ik}r_{j}, \\
  T^{S}_{ijk}(\mathbf{r}) &= \alpha^{2}\mathcal{T}_{1}^{L}(\alpha \|\mathbf{r}\|)(\delta_{jk}r_{i}+\delta_{ik}r_{j}+\delta_{ij}r_{k}) + \alpha^{4}\mathcal{T}_{2}^{L}(\alpha \|\mathbf{r}\|)r_{i}r_{j}r_{k}  + \alpha^{2}\mathcal{T}_{3}^{L}(\alpha \|\mathbf{r}\|)\delta_{ik}r_{j},\\
  T^{C}_{ij}(\mathbf{r}) &= \alpha^{2}r_{i}r_{j}/8\pi,\\
    T^{C}_{ij}(\mathbf{r}) &= -1/\pi.
\end{align}

\end{appendices}
\clearpage

\bibliographystyle{abbrvnat_mod}
\bibliography{largealpha}

\end{document}